\documentclass[lettersize,journal]{IEEEtran}
\usepackage{amsmath,amsfonts}
\usepackage{algorithmic}
\usepackage{algorithm}
\usepackage{array}
\usepackage{booktabs}
\usepackage[caption=false,font=normalsize,labelfont=sf,textfont=sf]{subfig}
\usepackage{textcomp}
\usepackage{stfloats}
\usepackage{url}
\usepackage{verbatim}
\usepackage{graphicx}
\usepackage[numbers, sort&compress]{natbib} 
\usepackage{times}
\usepackage{soul}
\usepackage{url}
\usepackage[hidelinks]{hyperref}
\usepackage[utf8]{inputenc}
\usepackage{caption}
\usepackage{adjustbox} % For resizing tables
\usepackage{makecell}
\usepackage{rotating}
\usepackage{graphicx}
\usepackage{amssymb}
\usepackage{amsthm}
\usepackage[switch]{lineno}
\usepackage{threeparttable}
\usepackage{multirow}
\usepackage{natbib}
\usepackage{rotating}
\usepackage{tikz,xcolor}
\usepackage{tabularx} % For tables that fit a specific width
\usepackage{lipsum}   % For dummy text to show the layout
\DeclareRobustCommand{\orcidicon}{
\begin{tikzpicture}
\draw[lime, fill=lime] (0,0)
circle[radius=0.16]
node[white]{{\fontfamily{qag}\selectfont \tiny \.{I}D}};
\end{tikzpicture}
\hspace{-2mm}
}
\foreach \x in {A, ..., Z}{%
\expandafter\xdef\csname orcid\x\endcsname{\noexpand\href{https://orcid.org/\csname orcidauthor\x\endcsname}{\noexpand\orcidicon}}
}
\begin{document}

\title{RideAgent: An LLM-Enhanced Optimization Framework for Automated Taxi Fleet Operations}

% \author{\IEEEauthorblockN{Anonymous Authors}}
\author{Xinyu Jiang, Haoyu Zhang, Mengyi Sha, Zihao Jiao, Long He, Junbo Zhang, Wei Qi
\thanks{
Xinyu Jiang and Mengyi Sha are with the Department of Industrial Engineering, Tsinghua University, Beijing, China (e-mail: jiangxin24@mails.tsinghua.edu.cn; mengyisha@mail.tsinghua.edu.cn).
Haoyu Zhang and Zihao Jiao are with the School of Computer and Artificial Intelligence, Beijing Technology and Business University, Beijing, China (e-mail: haoyuzhang@st.btbu.edu.cn; jiaozihao@btbu.edu.cn).
Long He is with the School of Business, George Washington University, Washington, District of Columbia (e-mail: longhe@gwu.edu).
Junbo Zhang is with JD Intelligent Cities Research (e-mail: msjunbozhang@outlook.com).
Wei Qi is with the Department of Industrial Engineering, Tsinghua University, Beijing, China, and also with the Desautels Faculty of Management, McGill University, Montreal, Quebec, Canada (e-mail: qiw@tsinghua.edu.cn).}}

% \markboth{IEEE TRANSACTIONS ON AUTOMATION SCIENCE AND ENGINEERING}%
% {Shell \MakeLowercase{\textit{et al.}}: RideAgent: How to Align Operations Research with Taxi Fleet Operators? An LLM-Enhanced Feature-Driven Optimization Framework}


\maketitle

\begin{abstract}
Efficient management of electric ride-hailing fleets, particularly pre-allocation and pricing during peak periods to balance spatio-temporal supply and demand, is crucial for urban traffic efficiency. However, practical challenges include unpredictable demand and translating diverse, qualitative managerial objectives from non-expert operators into tractable optimization models. This paper introduces RideAgent, an LLM-powered agent framework that automates and enhances electric ride-hailing fleet management. First, an LLM interprets natural language queries from fleet managers to formulate corresponding mathematical objective functions. These user-defined objectives are then optimized within a Mixed-Integer Programming (MIP) framework, subject to the constraint of maintaining high operational profit. The profit itself is a primary objective, estimated by an embedded Random Forest (RF) model leveraging exogenous features. To accelerate the solution of this MIP, a prompt-guided LLM analyzes a small sample of historical optimal decision data to guide a variable fixing strategy. Experiments on real-world data show that the LLM-generated objectives achieve an 86\% text similarity to standard formulations in a zero-shot setting. Following this, the LLM-guided variable fixing strategy reduces computation time by 53.15\% compared to solving the full MIP with only a 2.42\% average optimality gap. Moreover, this variable fixing strategy outperforms five cutting plane methods by 42.3\% time reduction with minimal compromise to solution quality. RideAgent offers a robust and adaptive automated framework for objective modeling and accelerated optimization. This framework empowers non-expert fleet managers to personalize operations and improve urban transportation system performance.

\end{abstract}
\def\abstractname{Note to Practitioners}
\begin{abstract}
The daily management of electric ride-hailing fleets, particularly the tasks of taxi pre-allocation and dynamic pricing during peak hours, poses operational complexities for service providers aiming to balance supply with urban demand. Translating diverse business objectives—such as maximizing revenue or enhancing taxi utilization—into actionable strategies can be challenging for managers without optimization expertise. This work introduces RideAgent, a Large Language Model (LLM)-based framework that automates and simplifies these core operational decisions. RideAgent assists fleet managers by interpreting their objectives stated in natural language into formal optimization goals. It then uses historical data to accelerate the computation of pre-allocation and pricing decisions. For fleet management practitioners, RideAgent offers a more intuitive and adaptive approach to operations: it enables the incorporation of specific managerial insights into the optimization process and significantly reduces the time needed to generate high-quality operational decisions. Our findings show that solution times can be reduced by over 50\% while maintaining near-optimal results. This allows for more agile and personalized decision-making. By lowering the technical barrier to using advanced optimization, RideAgent empowers a broader range of personnel to contribute to efficient fleet management. As the ride-hailing industry adopts AI-driven solutions, RideAgent provides a practical pathway for integrating advanced natural language understanding and automated optimization into daily fleet operations, improving urban mobility services.
\end{abstract}

\section{Introduction}
\IEEEPARstart{A}{s} electric taxi services become integral to urban mobility, fleet operators face the critical challenge of mitigating the spatio-temporal mismatch between taxi supply and passenger demand ~\cite{10.5555/3367032.3367057}. This requires the joint optimization of two interdependent strategies \cite{zhang2023}: taxi \textit{pre-allocation} to anticipate demand hotspots and dynamic \textit{pricing} to modulate demand.  
These strategies should be flexible enough to adapt to the ebb and flow of demand, yet robust enough to ensure profitability.

To mitigate these practical challenges, fleet operators increasingly turn to operations research (OR) techniques to jointly optimize zone-level pre-allocation and dynamic pricing. OR provides a set of rigorous methods for mathematical modeling, statistical inference, and algorithmic optimization for complex decision-making \cite{hillier2015introduction}. Classic methods—such as mixed-integer programming (MIP) and dynamic programming (DP)—can, in principle, deliver provably optimal resource allocation and pricing solutions. However, when these formulations must incorporate high-dimensional demand covariates, fine-grained spatial grids, and evolving operator objectives, they become exceedingly difficult to construct and computationally prohibitive to solve. Consequently, practitioners face the following limitations:

(i) \textbf{Data and Computational Challenges}. 
Accurate urban taxi demand forecasting and optimization are complex due to various external factors (e.g., weather, public events). Integrating these dynamic covariates into predictive models increases model complexity and computational overhead. Moreover, the trend towards finer-grained urban management, characterized by an expanding number of pre-allocation areas, imposes a substantial computational burden. Consequently, there is a pressing need for scalable algorithms that can solve large-scale optimization models in near real-time to ensure operational effectiveness.

(ii) \textbf{Dynamic and Personalized Objectives.} 
The paradigm of urban management is undergoing a fundamental shift from experience-based decision-making to a data-driven approach, propelled by the increasing availability of granular data within smart city ecosystems \cite{qi2019smart, hasija2020smart}. This transition has fostered a greater reliance on Key Performance Indicators (KPIs) and quantifiable objectives to guide operational strategies. Consequently, urban operators now demand optimization models that are not only contextually aware and actionable but also adaptable to a diverse set of objectives. These objectives often encompass competing priorities such as cost minimization, operational efficiency, service equity, and customer satisfaction. Traditional OR models, which typically feature fixed and predetermined objective functions, often lack the flexibility to accommodate these dynamic and personalized preferences.

(iii) \textbf{Knowledge Gaps.}   
A significant barrier to the widespread adoption of advanced OR techniques in urban mobility lies in the knowledge gap between OR specialists and operational practitioners. Complex optimization models pose a steep learning curve for non-expert operators, limiting their effective leveraging of OR tools. Conversely, this gap also complicates the task for OR modelers, who struggle to accurately formulate objective functions that are nonlinear, nonconvex, and truly representative of the practitioners' real-world goals.

In light of these obstacles, frameworks that couple \emph{Large Language Models} (LLMs) with domain-specific optimization engines offer a promising way forward by leveraging their capabilities in natural language understanding, reasoning, and tool use\,\cite{shen2024hugginggpt}.
However, most current LLM–OR systems still rely on shallow natural-language parsing followed by a one-shot call to a generic solver\,\cite{ijcai2024p212}; this often compromises modeling fidelity and yields biased or even infeasible solutions.
For instance, even after prompt tuning, GPT-4 reports an 11\,\% optimality gap on a 50-node traveling-salesperson instance\,\cite{yang2023large,achiam2023gpt}.
Moreover, the enlarged, tightly coupled models required for real-time joint pre-allocation and pricing remain difficult to solve within operational deadlines. 

To address these limitations, we introduce \textbf{RideAgent}, an LLM-powered framework designed to make large-scale fleet optimization both accessible and computationally tractable. RideAgent employs an LLM in two capacities: problem formulation and solution acceleration. First, it translates the natural language objectives of non-expert operators into precise mathematical formulations for the optimization model. Our primary innovation, however, is a novel LLM-guided acceleration heuristic that performs model reduction. Distinct from approaches that use LLMs to construct models from scratch, we start with a full-scale, expert-defined optimization model and use the LLM to intelligently prune it. By learning from few-shot examples of historical optimal solutions, the agent identifies and fixes decision variables with low sensitivity to the optimal solution. This intelligent reduction of the decision space renders the complex, high-dimensional optimization problem more tractable for conventional solvers like Gurobi. The result is a hybrid architecture that dramatically accelerates computation for real-time decision-making while maintaining near-optimal solution quality.

\begin{itemize}
    \item \emph{An LLM-Powered Architecture for Feature-Driven Fleet Optimization}. In order to account for environmental uncertainties and covariates, RideAgent incorporates a feature-driven optimization framework that integrates a random forest (RF) with an optimization model. This framework allows for the inclusion of numerous features, thereby enhancing the alignment with user requirements and reducing decision-making bias. 
    \item \emph{An LLM-Guided Variable Fixing Heuristic for Accelerated Optimization}. We introduce a novel variable fixing strategy guided by LLMs' reasoning ability. Based on limited historical optimal decision data, RideAgent can learn the rules of optimal decision-making and proactively fix a subset of decision variables that have low sensitivity to the optimal solution. By focusing on key variables, this approach simplifies the pre-defined large-scale optimization problem and enhances the computational efficiency.
    \item \emph{Real-World Case Study Validation of LLM-Driven Optimization}. We validate the effectiveness of RideAgent through a comprehensive real-world case study. Results demonstrate that RideAgent excels in generating accurate objective functions to accommodate user requirements. 
    Benchmarking against the full-scale MIP model and established cutting plane methods, RideAgent demonstrates significant acceleration with minimal compromise in solution quality. 
\end{itemize}
The remainder of this paper is organized as follows: Section 2 reviews some related work. Section 3 defines the optimization models. Section 4 introduces the structure and functions of RideAgent. Section 5 analyzes the results of the case study. Section 6 summarizes our work.

\section{Related Work}
Mitigating the supply-demand mismatch in urban taxi systems is a cornerstone of enhanced urban mobility. A substantial body of literature has focused on this challenge by optimizing joint taxi pre-allocation and pricing strategies. However, the integration of LLM agents to simultaneously address modeling flexibility and computational scalability in this domain remains unexplored. This section reviews the pertinent literature across three key areas to situate our contributions.

\subsection{Joint Taxi Pre-Allocation and Pricing}
The joint taxi pre-allocation and pricing problem has attracted considerable attention in the transportation literature~\cite{fanti2021innovative, liu2025}. 
The intrinsic fluctuations of service demand introduce complexity to the decision-making framework ~\cite{ozkan2020joint,its2020}. Conventional deterministic optimization models do not account for dynamic changes, which may result in biased solutions. To address these challenges, an increasing number of researchers have adopted various optimization approaches to account for spatiotemporal uncertainty in service demand and enhance decision-making effectiveness. Existing studies propose stochastic or robust optimization models ~\cite{xu2018electric, pantuso2022exact, relocation2022tase, EJOR2024} or dynamic programming models ~\cite{shah2022joint,chen2023real}, and then use reformulation techniques, heuristic algorithms or approximation algorithms to derive efficient solutions. 

Existing research often overlooks the impact of covariate features on uncertainty and struggles to incorporate features into OR models \cite{hao2020robust}. Additionally, the complexity of the joint strategy necessitates OR expertise, which may limit its wider application to non-expert groups. To bridge these gaps, we integrate OR models with LLMs to provide more accessible and reliable solutions.

\subsection{LLM-Assisted OR Tools}
LLMs are increasingly being explored for their potential to augment OR. Current research predominantly follows two primary approaches:

(1) \textbf{LLMs as Optimizers}. This paradigm utilizes techniques like in-context learning \cite{nie2024importance, ORLM2025} and Chain-of-Thought (CoT) prompting \cite{kojima2022large} to solve optimization problems directly through conversational interfaces. A key appeal of this approach is its accessibility, as it often does not require users to have specialized OR expertise.
%This technique leverages the in-context learning LLM (Dong et al. 2022) to tackle optimization challenges directly via conversational interfaces. By integrating chain of thought (CoT) prompts (Kojima etal. 2022), the technique can utilize LLM’s reasoning abilities to solve the optimization problems and only needs a few prompts (Wei et al. 2022) related to the users’ problems. In this way, LLMs as optimizers do not require extensive knowledge of specific details or specialized OR expertise.

(2) \textbf{LLM Agents} \cite{zhao2024expel, huang2024large}.
This approach integrates LLMs with external tools, such as callable APIs \cite{chacon2024large, xu2024large} or traditional OR solvers, to form a more robust system. By offloading complex calculations to specialized tools, these agents can mitigate the hallucination and token-limitation issues inherent in standalone LLMs. Recent studies in this domain have focused on creating agents that can solve MIPs \cite{ahmaditeshnizi2023optimus, orllm2025, zhang2025decision}.

Nevertheless, these emerging techniques still face significant limitations in terms of solution quality and modeling flexibility. LLMs as optimizers, for instance, often struggle to guarantee the accuracy or optimality of solutions for complex combinatorial problems \cite{yang2023large}. Meanwhile, existing LLM agents, especially those interfacing with MIP solvers \cite{ahmaditeshnizi2023optimus, li2023synthesizing}, are often constrained to static and well-defined problems. They typically lack the flexibility to dynamically formulate the nuanced or non-standard objective functions required to address the diverse and evolving requests of practitioners.

\subsection{Feature-Driven Optimization}
While LLM-powered agents can formulate optimization models, they often rely on deterministic frameworks (e.g., standard MIPs). Such models struggle to account for the inherent uncertainties of dynamic environments, like urban mobility, potentially leading to biased or suboptimal solutions when faced with fluctuating demand or changing external conditions.
The Feature-Driven Optimization paradigm directly addresses this challenge. Instead of solving a problem based on average or point estimates of uncertain parameters, feature-driven optimization seeks to learn a policy or decision rule that maps observable covariates (features) directly to optimal decisions \cite{review2024}. This approach integrates machine learning with optimization, allowing decisions to adapt to real-time information. Dominant feature-driven paradigms include decision rule optimization, sequential learning and optimization, and integrated learning and optimization ~\cite{qi2022integrating,review2024}.
Our proposed RideAgent adopts the latter approach and extends the integrated framework of Biggs et al. \cite{Biggs_2022}. Specifically, we embed predictions from a pre-trained RF directly into the MIP's objective function. This allows our model to leverage a rich set of exogenous features (e.g., weather, time of day, public events) and makes the core profit-maximization objective highly responsive to real-world conditions. The result is a more robust optimization engine that yields contextually relevant and practically grounded pre-allocation and pricing strategies.

\section{Problem Formulation}\label{Problem Formulation}
\subsection{Joint Taxi Pre-Allocation and Pricing Model}
We address the optimal pre-positioning of an electric taxi fleet to mitigate spatio-temporal supply-demand mismatches. The mismatch results in service shortages in some areas and idle taxis in others. To formalize this, we first partition the city into a set of operational areas. Areas with a surplus of taxis relative to local demand are defined as the set of \textbf{supply areas}, $I = \{1, \dots, |I|\}$. Conversely, areas with a service deficit form the set of \textbf{demand areas}, $J = \{1, \dots, |J|\}$. Our model aims to mitigate this supply-demand imbalance by optimizing taxi repositioning from set $I$ to set $J$. Electric taxis are categorized by their state of charge (SOC), which is discretized into a set of levels $K = \{1, \dots, |K|\}$, ordered from lowest to highest. A key operational rule is that a request for a taxi with SOC $k$ can be served by any available taxi with an SOC of $k$ or higher, but not lower.

Our model setup is informed by the framework for taxi pre-allocation problems presented in Hao et al. (2020) \cite{hao2020robust}. The system is defined by a set of given parameters and the key decisions the operator must make. The primary inputs include the initial supply of taxis $S_{ik}$ in each area and the parameters that affect the demand. The pre-allocation cost of each taxi paid by operators is $w_{ij}=\hat{w}_{ij}+b_j$, where $\hat{w}_{ij}$ is the inconvenience cost and $b_j$ is the fixed online booking fee. $\hat{w}_{ij}$ is proportional to the distance between areas. The total fare paid by a customer is composed of a variable mileage-based fee $\hat{u}_{jk}$ and the booking fee $b_j$. For each successful trip, the booking fee functions as a passthrough payment: the operator collects $b_j$ from the customer and transfers it directly to the driver who is allocated to demand area $j$. Operators earn average revenue $u_{jk}$ per order for taxis with SOC $k$: $u_{jk}=\theta \cdot \hat{u}_{jk}+b_{j}$. 

Given these parameters, the operator's core task involves two sets of decisions: first, determining the number of taxis $x_{ijk}$ to reposition, and second, setting the variable portion of the fare $\hat{u}_{jk}$. 
The operating income is the order revenue in demand areas: $R(\mathbf{\hat{u}},\mathbf{d})=\sum_{j \in J}\sum_{k \in K}u_{jk} d_{jk}$. The operating cost is the taxi pre-allocation cost: $C(\mathbf{x})=\sum_{i \in I} \sum_{j \in J}\sum_{k \in K} w_{ij} x_{ijk}$. The goal is to maximize the operational profit: $R(\mathbf{\hat{u}},\mathbf{d})-C(\mathbf{x})$.
To model the cascading nature of SOC fulfillment, we introduce an auxiliary variable $v_{jk}$ that represents the number of taxis with an SOC greater than $k$ that are available in area $j$ after demands for all higher SOC levels have been met. All symbols and their definitions are consolidated in Table~\ref{tab:notation}.

% \captionsetup[table]{name={TABLE},labelsep=space}
\begin{table}[htbp]
\centering
\caption{Notation for the Problem Formulation}
\label{tab:notation}
\begin{tabularx}{\columnwidth}{@{}l X@{}}
\toprule
\textbf{Symbol} & \textbf{Definition} \\
\midrule
$I, J, K$ & Sets of supply areas, demand areas, and SOC levels, respectively. \\
$S_{ik}$ & Number of available taxis with SOC $k$ at supply area $i$. \\
$z_{jk}$ & Anticipated demand for SOC $k$ or higher in demand area $j$. \\
$w_{ij}$ & Unit cost of pre-allocating a taxi from area $i$ to area $j$. \\
$b_j$ & Fixed booking fee for a trip originating in area $j$. \\
$\theta$ & Operator's revenue share from the variable fare. \\
\midrule
&\textbf{Decision Variables} \\
$x_{ijk}$ & Number of taxis with SOC $k$ pre-allocated from $i$ to $j$. \\
$\hat{u}_{jk}$ & Variable portion of the average fare for a trip from $j$ with SOC $k$. \\
\midrule
&\textbf{Auxiliary Variables} \\
$d_{jk}$ & Number of satisfied demands for SOC $k$ in demand area $j$. \\
$v_{jk}$ & Surplus taxis from higher SOC levels available for level $k$ in area $j$. \\
\bottomrule
\end{tabularx}
\end{table}

Based on these definitions, we formulate the electric taxi pre-allocation and pricing optimization model as follows. For clarity, we use the notation $[A]^+ = \max(0, A)$ and $A \wedge B = \min(A, B)$.
\begin{align}
\max\limits_{\mathbf{x},\mathbf{\hat{u}},\mathbf{v},\mathbf{d}} \quad & R(\mathbf{\hat{u}},\mathbf{d})-C(\mathbf{x}) \label{obj-origin}\\
\text{s.t.} \quad & \sum_{j \in J}x_{ijk} \le S_{ik},\forall i\in I, k\in K, \label{con1}\\
&d_{jk}=z_{jk}\wedge \big[\sum_{i \in I}x_{ijk}+v_{j(k+1)}\big],\forall j\in J, k\in K,\label{con2}\\
&v_{jk}=\big[\sum_{i \in I} x_{ijk}-z_{jk}+v_{j(k+1)}\big]^+ ,\forall j\in J, k\in K,\label{con3}\\
&v_{j(|K|+1)}=0,\forall j\in J, \label{con4}\\
&x_{ijk} \in \mathbb{N}, \forall i\in I, j\in J, k\in K. \label{con5}
\end{align}

Objective (\ref{obj-origin}) represents the total operational profit from taxi pre-allocation and pricing decisions. Constraint (\ref{con1}) ensures that the number of taxis allocated from a supply area does not exceed its available fleet. Constraint (\ref{con2}) then determines the satisfied demand based on the total available supply, which is composed of the newly allocated taxis and the surplus taxis from higher SOC levels. Constraint (\ref{con3}) calculates the resulting surplus for the next level. Constraint (\ref{con4}) sets the boundary condition for this process, while Constraint (\ref{con5}) enforces the integrality of the allocation decisions.
\subsection{RF-based Feature-Driven Model}
The model formulated in the previous section operates under a predict-then-optimize paradigm. Demand is predicted upfront and then treated as a fixed input for the optimization stage. A well-known limitation of such sequential frameworks is error propagation: any inaccuracies from the initial demand prediction are inevitably carried into the optimization model, which may then yield a solution that is optimal for the flawed prediction but suboptimal in reality. To circumvent error propagation, we replace the deterministic objective function with a feature-driven profit prediction model based on a RF. This model directly learns the relationship from system features and operator decisions to the final predicted profit, thereby avoiding the intermediate prediction step and its associated error accumulation. To embed this new objective model into our MIP, we apply the RF-to-MIP conversion technique from Biggs et al. (2022) \cite{Biggs_2022}.

% Deterministic models are generally ineffective because external factors have a significant impact on user demands.
% Therefore, we introduce a RF-based feature-driven model, which can make full use of historical operation data and explore the relationship between system features (such as weather and operation decisions) and operation results (profit). The tree structure also augments LLM's reasoning process. In addition, the feature-driven structure can achieve direct prediction from decisions to targets, reducing the accumulation of random parameter prediction errors.

\begin{figure}[htbp]
\centering
\includegraphics[width=\columnwidth]{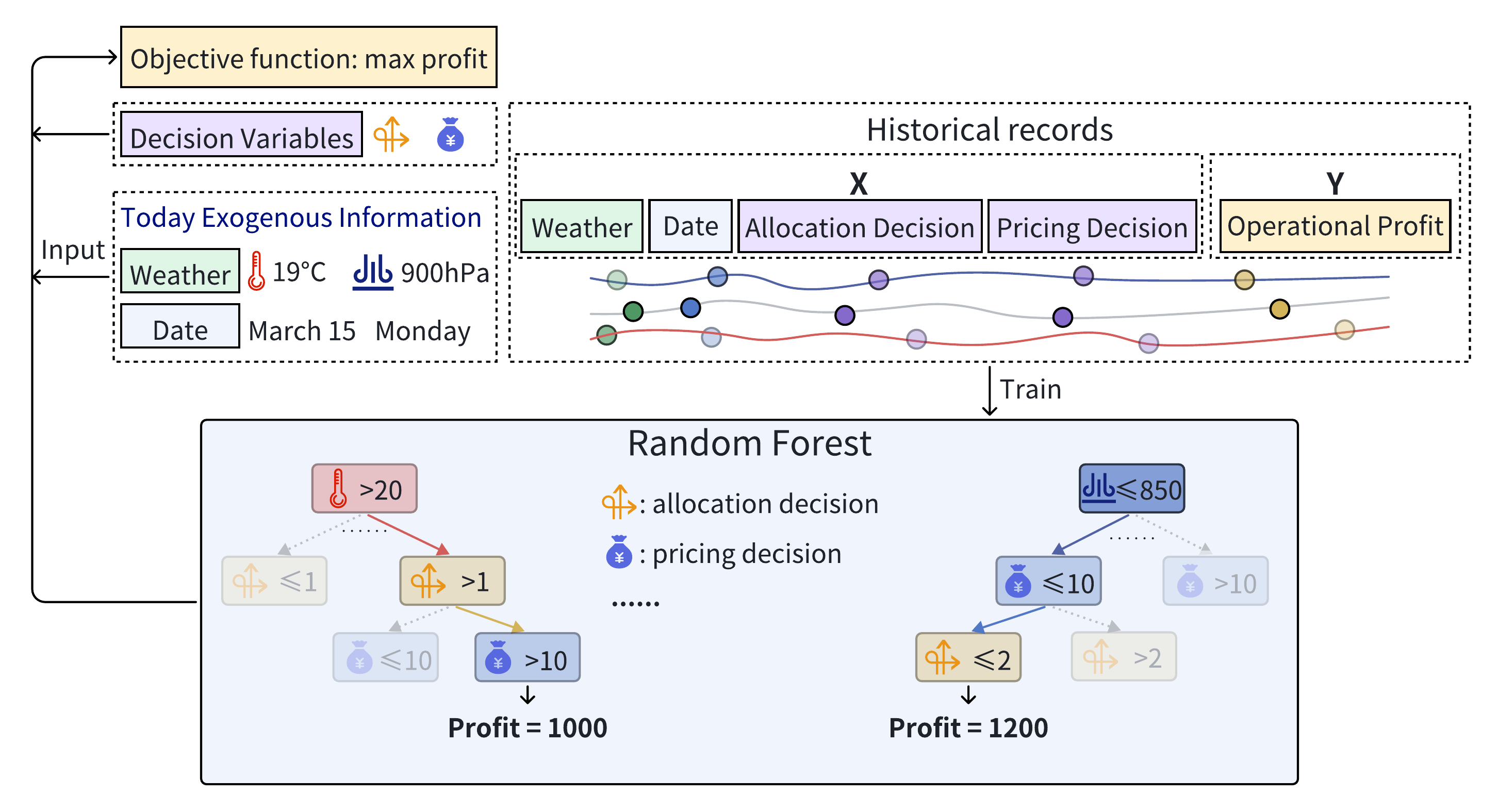}
\captionsetup{name={Fig.},labelsep=period,justification=raggedright, singlelinecheck=false}
\caption{Training and Usage Process of the random forest}
\label{fig:rf}
\end{figure}

The core of this integrated method is to embed the entire trained RF structure as a set of constraints within a MIP. Figure~\ref{fig:rf} illustrates the workflow of training and deploying the RF model. The RF is first trained using historical records, which include exogenous features like weather and date, and operational data such as allocation and pricing decisions, along with the resulting profit. In the second stage, the trained model is deployed for optimization. It takes the current day's exogenous information as fixed inputs. The allocation and pricing actions are treated as decision variables, and an optimization solver selects their optimal values. This process effectively navigates a path through the tree's branches to a leaf node that predicts a final profit. The objective is to determine the set of decisions that leads to the leaf with the maximum profit.

% We illustrate a sample tree of RF in Figure \ref{Figure RF}. 
Suppose there are $|H|$ trees in the RF, and the leaf node of each tree $h$ represents the predicted operational profit $P(\mathbf{y},\mathbf{c})$, where $\mathbf{y}=(\mathbf{x}, \mathbf{\hat{u}})$ is the decision variable and $\mathbf{c}$ is the feature. Our goal is to maximize the operation profit by maximizing $P(\mathbf{y},\mathbf{c})$. Let $N^{h}$ denote the number of nodes (excluding the leaves) in tree $h$. For each interior node $n$, let $p_n, l_n$ and $r_n$ be the immediate parent, the left and right children, respectively. Let $L^h$ represent the set of leaves in tree $h$, and let $P_m^{h}\;(m\in L^h)$ denote the score of each leaf. Given the trained RF, we introduce binary variables $q_{n_1,n_2}^{h}$ to select branches and decide the range of the variables. The resulting feature-driven optimization model is as follows:
\begin{align}
\max \limits_{\mathbf{y},\mathbf{q}}\quad & \frac{1}{|H|}\sum_{h=1}^{|H|}\sum_{m\in L^h}P_m^{h}q_{p_m,m}^{h} \label{obj-RF}\\
\text{s.t.}\quad & a_{n, h} \mathbf{y}-M\left(1-q_{n, l_n}^h\right) \leq b_n^h, \forall h \in H, n \in N^h,\label{con21} \\
& a_{n, h} \mathbf{y}+M\left(1-q_{n, r_n}^h\right) \geq b_n^h, \forall h \in H, n \in N^h,\label{con22}\\
& q_{n, l_n}^h+q_{n, r_n}^h=q_{p_n, n}^h, \forall h \in H, n \in N^h,\label{con23} \\
& \sum_{n \in L^h} q_{p_n, n}^h=1, \forall h \in H,\label{con24}\\
& q_{n, l_n}^h, q_{n, r_n}^h, q_{p_n, n}^h \in\{0,1\}, \forall h \in H, n \in N^h, \label{con25} \\
& \text{Constraints} \; (\ref{con1})\text{-}(\ref{con5})\nonumber.
\end{align}

Objective (\ref{obj-RF}) is to maximize the average predicted profit across all trees in the forest. Constraint (\ref{con21}) and 
constraint (\ref{con22}) are big-M logical constraints to determine which leaf the solution $\mathbf{y}$ lies in. Constraint (\ref{con23}) ensures that if a parent node is inactive, its children must also be inactive; but if any child is active, then the parent must also be active. Constraint (\ref{con24}) guarantees that within each tree $h$, only one leaf can be active. Constraint (\ref{con25}) defines the binary variables. Finally, the solution must also satisfy the operational constraints (\ref{con1})\text{-}(\ref{con5}) from the previous section to ensure its physical feasibility.
% \begin{figure}[htbp]
%     \begin{minipage}[b]{0.5\textwidth}
%       \centering
%       \includegraphics[width=0.95\textwidth]{Figure/tree1.png}\vspace{0.1cm}
%                         \begin{flushleft}
%       \small \centering (a) Sample tree of the trained random forest.
%             \end{flushleft}
%     \end{minipage}%
%     \mbox{\hspace{-0.2cm}}
%     \begin{minipage}[b]{0.5\textwidth}
%       \centering
%       \includegraphics[width=0.95\textwidth]{Figure/tree2.png}\vspace{0.1cm}
%                              \begin{flushleft}
%       \small \centering (b) Sample tree after variable fixing.
%             \end{flushleft}
%     \end{minipage}
%     \caption{Sample Tree of the Trained Random Forest for before and after variable fixing.}
%     \label{fig_RF}
% \end{figure}
\subsection{LLM-embedded Feature-Driven Model}
% The LLM can use its reasoning ability to understand user's ambiguous questions and add a new objective function for the basic MIP model. 
% In addition, the LLM can infer the decision variables that contribute less to the optimal solution from the historical optimal strategy with fewer samples. Parameterizing the decision variables selected by the LLM as the mean of the historical optimal decision can effectively shorten the solution time without sacrificing too much accuracy of the solution.
% After receiving and analyzing the user's question, the LLM agent will convert the user's question text into the target function code that can be run in the solver.
While our feature-driven model can find profit-maximizing strategies, it faces two key practical challenges in large-scale urban operations. First, operators often have complex and dynamic goals that are difficult to translate into a mathematical objective. Second, the sheer scale of urban optimization problems often makes finding exact solutions computationally prohibitive in a reasonable timeframe. To address both challenges, we introduce an LLM-powered agent that performs two synergistic functions: 1) dynamic objective formulation to translate user queries into formal objectives, and 2) heuristic model reduction to accelerate computation.

To accelerate computation, the agent employs a novel heuristic strategy learned from a small set of \textbf{pre-solved optimal instances}. By analyzing these few-shot examples, which pair historical scenarios with their true optimal decisions, the LLM learns to identify a subset of decision variables $\mathbf{y}^{\prime} \subset \mathbf{y}$ that consistently exhibit low sensitivity to shaping the optimal solution. Then it proposes fixing these variables to their historical average values, $\bar{\mathbf{y}}^{\prime}$, leaving a smaller set of active decision variables, $\hat{\mathbf{y}} = \mathbf{y} \setminus \mathbf{y}'$, to be optimized. This model refinement process is achieved by employing a prompted LLM, which 
we formally denote as $\textbf{LLM}(Q;\mathcal{D},P^{\text{PP}})$. In this formulation, a structured prompt $P^{\text{PP}}$, which incorporates the user's query $Q$, instructs the LLM to analyze historical data 
$\mathcal{D}$ as the evidence base for identifying low-sensitivity variables~\cite{salemi2024optimization}. The LLM agent is described in detail in the next section. %Rather than optimizing the operations of the overall city, RideAgent optimizes the parameterized model. In summary, RideAgent can leverage the logical reasoning ability of LLMs to generate the objective function that users care about, and learn from the historical optimal decisions of smaller samples to reduce the scale of the problem.
This dual-functionality results in a bi-objective optimization problem that we solve using a \textbf{lexicographical method}. The primary goal is to maximize the predicted profit. The secondary goal is to optimize the new LLM-generated objective, $\textbf{LLM}(Q;\mathcal{D},P^{\text{IG}})$, where $P^{\text{IG}}$ denotes the predefined prompts.
Define $A(\mathbf{c})F(\mathbf{y})\leq K(\mathbf{c})$ to depict the feasible region and constraints related to $\mathbf{y}$. After reducing the scope of the RF-based model by fixing some decision variables, we have the following optimization problem:
\begin{align}
\max\limits_{\hat{\mathbf{y}},\mathbf{q}} \quad & \frac{1}{|H|}\sum_{h=1}^{|H|}\sum_{m\in L^h}P_m^{h}q_{p_m,m}^{h} \label{obj-rf}\\
\min_{\hat{\mathbf{y}},\mathbf{q}}\quad &\textbf{LLM}(Q;\mathcal{D},P^{\text{IG}}) \label{obj-llm}\\
\text{s.t.}\quad & A(\mathbf{c})F(\hat{\mathbf{y}}, \bar{\mathbf{y}}^{\prime})\leq K(\mathbf{c}), \label{con31}\\
& \text{Constraints} \; (\ref{con21})\text{-}(\ref{con25})\nonumber.
\end{align}

Objective (\ref{obj-llm}) is generated by the LLM agent based on the user's query $Q$ and is treated as a secondary objective, optimized only after the primary, RF-based objective (\ref{obj-RF}) is maximized. %In the multi-objective optimization problem, we solve the programming model with objective (\ref{obj-rf}) first, and then re-optimize the problem with the newly added objective (\ref{obj-llm}) while maintaining the first objective optimal.
Constraint (\ref{con31}) compactly represents the feasible region of the problem, where $\hat{\mathbf{y}}$ and the fixed variables $\bar{\mathbf{y}}^{\prime}$ must satisfy all operational constraints. Constraints (\ref{con21})\text{-}(\ref{con25}) involve binary variables to model the branch and leaf selection in RF.

After fixing some redundant variables $\mathbf{y}^{\prime}$, RideAgent offers a more agile model.% , especially notable for its reduced computation time. % We will evaluate the performance by examining both the optimality gaps and the time efficiency of the computations in our case study.

% \begin{figure}[htbp]
% \centering
% \includegraphics[width=\columnwidth]{Figure/tree2.png}
% \caption{A sample tree after variable fixing}
% \label{fig_RF}
% \end{figure}
\section{LLM-based Agent}
% This section introduces the structure of our proposed agent and its various functions. The macro structure of the agent adopts a process similar to human cooperation and decision-making, and adds some tricks in specific functions to adapt to the taxi pre-allocation and pricing problem ~\cite{wei2022chain,xiao2023chain}. A visual representation of the agent framework is provided in Figure \ref{Figure agent}.  %Please refer to Appendix \ref{app:pse} for the detailed pseudocode of the agent implementation process.

% This section first presents the architecture and workflow of our proposed \textit{RideAgent}. We then elaborate on the core intuition behind its design, which synergizes the reasoning capabilities of LLMs with formal optimization methods to solve complex urban operations problems.
This section details the architecture and workflow of our proposed \textit{RideAgent}, the framework responsible for the model's agility and for orchestrating the overall optimization process. We first present its multi-component framework, which synergizes the reasoning capabilities of LLMs with formal optimization methods  We then conclude by discussing the core intuition behind its design, which we term Small-Sample Guided Optimization. 

\begingroup
\begingroup
\setlength\abovedisplayskip{0pt}
\setlength\belowdisplayskip{0pt}
\setlength\abovecaptionskip{5pt}
\setlength\belowcaptionskip{0pt}
\setlength\abovedisplayshortskip{0pt}
\setlength\belowdisplayshortskip{0pt}
\begin{figure*}
\centering
\includegraphics[width=\textwidth]{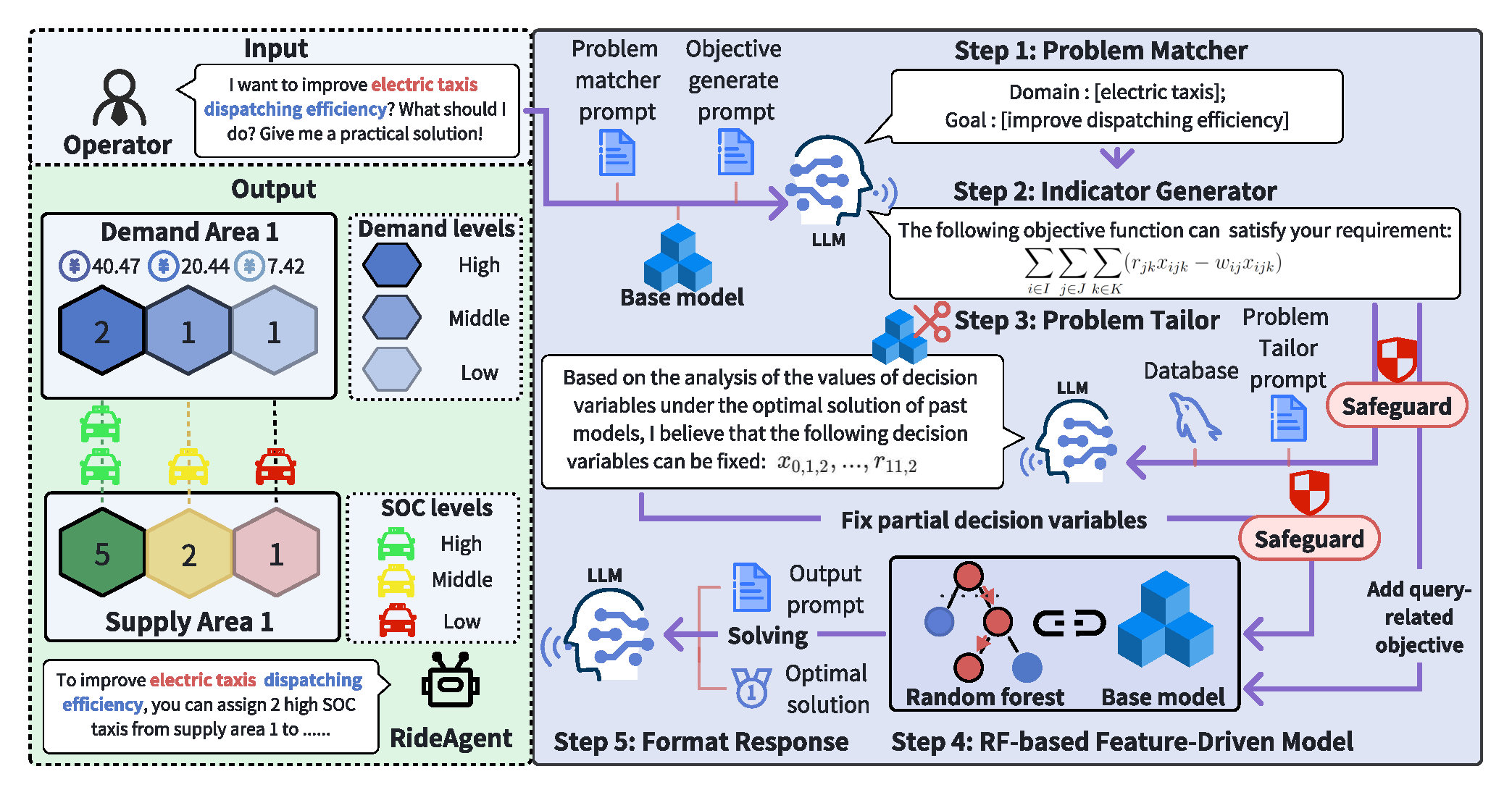} 
\captionsetup{name={Fig.},labelsep=period,justification=raggedright, singlelinecheck=false}
\caption{The agent framework}
\label{Figure agent}
\end{figure*}
\endgroup
\endgroup

\subsection{The RideAgent Framework}
The macro-structure of RideAgent mimics a human-like cooperative decision-making process, tailored for the taxi pre-allocation and pricing problem ~\cite{wei2022chain,xiao2023chain}. The agent's architecture, visualized in Figure \ref{Figure agent}, is composed of five key components that sequentially process a user's natural language (NL) query to produce an optimized operational decision. The detailed workflow is formally delineated in Algorithm~\ref{pse:1}.

\textbf{Step 1: Problem Matcher}. The process begins when a user submits a natural language (NL) query $Q$ focused on a particular objective within a defined region (e.g., \textit{How to improve the electric taxi dispatching efficiency?}). Such a query is then analyzed by Problem Matcher (see Figure \ref{Figure agent} \textit{Step 1}) to identify an appropriate domain-specific agent (e.g., \textit{electric taxis}) and routes the request for further processing. 

\textbf{Step 2: Indicator Generator}. Using a predefined prompt $P^{\text{IG}}$, the indicator generator (see Figure \ref{Figure agent} \textit{Step 2}) formulates NL queries $Q$ into a mathematical function $f(\mathbf{y};\mathbf{w})$ based on operation decisions $\mathbf{y}$ (e.g., number of allocated electric taxis and average order price for high SOC taxis) and parameters $\mathbf{w}$ in the feature-driven programming. 
\begin{align}
    f(\mathbf{y};\mathbf{w}) = \textbf{LLM}(Q;\mathcal{D},P^{\text{IG}}): \mathbb{R}^{|\mathcal{Y}|\times |\mathcal{W}|}\rightarrow \mathbb{R}_+,
\end{align}
where $\mathcal{Y},\mathcal{W}$ denote the feasible domain of decision variables $\mathbf{y}$ and set of parameters $\mathbf{w}$, respectively.
The objective function $f(\mathbf{y};\mathbf{w})$ is then integrated into the objective of MIP for re-optimization purposes.
Furthermore, the historical operational database $\mathcal{D}$ is available in the form of an API that generates Structured Query Language (SQL) by an LLM. $\mathcal{D}$ contains rich historical taxi operational data and a small amount of historical optimal decisions that maximize operating profits. 
Finally, \textbf{Code Safeguard} scrutinizes the output code of $f(\mathbf{y};\mathbf{w})$ to ensure compliance with the coding standards required by solvers like Gurobi and CPLEX \cite{he2024can}.

\textbf{Step 3: Problem Tailor}. 
The role of the Problem Tailor is to prune the vast decision space by identifying the most influential decision variables relevant to the query $Q$. A metric $S$ is established to evaluate how well a given solution satisfies the query $Q$. This score is typically defined as the percentage improvement of the indicator function $f$ for a new solution $\mathbf{y}^*_t$ compared to a baseline, where all decision variables are set to their historical average values $\bar{\mathbf{y}}_{hist}$:
$S_t=\frac{f(\hat{\mathbf{y}}^*_t,\bar{\mathbf{y}}^{\prime};\mathbf{w})-f(\bar{\mathbf{y}}_{hist};\mathbf{w})}{f(\bar{\mathbf{y}}_{hist};\mathbf{w})}$.

This module operates inside the iterative loop. The prompt $P_t^{\text{PP}}(\hat{\mathbf{y}}_{t-1},S_{t-1})$ in iteration $t$ contains previously remaining variables $\hat{\mathbf{y}}_{t-1}$ and satisfaction score $S_{t-1}$. Based on this information, the LLM provides guidance on which subset of variables $\mathbf{y}'_t$ is least sensitive and can be fixed to their historical average values $\bar{\mathbf{y}}'_t$. This leaves a smaller, more promising set of active variables $\hat{\mathbf{y}}_t = \mathbf{y} \setminus \mathbf{y}'_t$ for the solver to focus on.
%Instead of optimizing over all possible variables, it intelligently selects a smaller subset $\hat{\mathbf{y}}$.
% The role of Problem Tailor is to distill the broad queries $Q$ into more focused decisions subsets $\hat{\mathbf{y}}$ and to identify the pertinent factors. %, instead of examining the entire set of decisions variables within the feature-driven model's optimization problem $\min_{\mathbf{y}\in \mathcal{Y}} g(\mathbf{y};\mathbf{w})$.
% Problem Tailor initiates by utilizing SQL commands generated by the LLM to extract relevant information from the database $\mathcal{D}$, which contains rich historical taxi operational data and a small amount of historical optimal decisions that maximize operating profits. It then assesses which operational decisions contribute the least to the optimal solutions and subsequently retains the most relevant decision variables $\hat{\mathbf{y}}$. 
This interaction can be described by the following functional mapping:
\begin{align}
    \hat{\mathbf{y}}_t=\textbf{LLM}(Q;\mathcal{D},P_t^{\text{PP}}(\hat{\mathbf{y}}_{t-1},S_{t-1})).
\end{align}
%=\frac{|(f(\hat{\mathbf{y}}^*_{t-1},\mathbf{y}^{\prime};\mathbf{w})-f(\hat{\mathbf{y}};\mathbf{w})|}{f(\hat{\mathbf{y}};\mathbf{w})}
% This framework enables the LLM to evaluate the most pertinent decision variables $\hat{\mathbf{y}}$ and strive for an enhanced $S_{t}$ value.
% The remaining, $\mathbf{y}^{\prime}$ is aligned with average historical values $\mathbf{y}^\text{hist}$ and is represented as equivalent constraints in subsequent steps. 
After that, \textbf{Code Safeguard} reviews the constraints generated by the LLM to ensure they adhere to the syntax requirements of optimization solvers.

\textbf{Step 4: RF-based Feature-Driven Model}. %The RF-based feature-driven model programming encompasses decision variables pertinent to a specific problem along with associated features and parameters, collectively termed as covariates. 
The RF-based feature-driven model $\max_{\mathbf{y}\in \mathcal{Y}}g(\mathbf{y};\mathbf{w})$ (see details in the previous section) consists of two components (Figure \ref{Figure agent} \textit{Step 4}): (i) Pre-trained RF, tailored to a particular urban optimization goal (taxi operation profits), functions as a predictive analytics tool.
(ii) MIP: Drawing on previous research ~\cite{Biggs_2022}, the structure of the RF is embedded within an MIP model, as formulated in Section \ref{Problem Formulation}. The MIP's objectives encompass operational targets as well as indicators related to the query $f(\mathbf{y};\mathbf{w})$ supplied by Indicator Generator.
% The optimization of new objective $f(\mathbf{y};\mathbf{w})$ is addressed with secondary priority, ensuring that the original MIP objective is optimized first. 
% Furthermore, the fixed variables constraints$\mathbf{y}^{\prime}$, presented as equality constraints, are integrated into the MIP constraints.
The MIP takes the reduced set of active variables $\hat{\mathbf{y}}_t$ and the fixed variable constraints ($\mathbf{y}'_t = \bar{\mathbf{y}}'_t$) and solves the bi-objective optimization problem lexicographically. Subsequently, the agent proceeds to optimize the problem on a more focused decision subset: 
\begin{align}
\hat{\mathbf{y}}^*_t=\arg \max_{\hat{\mathbf{y}}} [g(\hat{\mathbf{y}},\bar{\mathbf{y}}_t^{\prime};\mathbf{w}), f(\hat{\mathbf{y}},\bar{\mathbf{y}}_t^{\prime};\mathbf{w})].
\end{align}
%By fixing irrelevant decisions $\mathbf{y}^{\prime}$, the reduced scope model $\max_{\hat{\mathbf{y}}} g(\hat{\mathbf{y}},\bar{\mathbf{y}}^{\prime};\mathbf{w})$ efficiently yields high quality solutions $\hat{\mathbf{y}}^*_t$ within reduced computational time.

\textbf{Step 5: Response Prompter}.  
The framework iteratively performs steps 3 and 4, generating a sequence of improving solutions $(\hat{\mathbf{y}}^*_t, S_t)$. After the loop terminates (i.e., when the satisfaction score $S_t$ no longer improves), Response Prompter translates the optimal solutions and associated scores $(\mathbf{y}^*,S^*)$ into a clear, understandable natural language summary for the end users.

\subsection{Core Intuition: Small-Sample Guided Large-Scale Optimization}

\begin{algorithm}[htbp]
        \caption{RideAgent pseudocode} 
        \renewcommand{\algorithmicrequire}{\textbf{Input:}}
        \renewcommand{\algorithmicensure}{\textbf{Output:}}
        \begin{algorithmic}[1]
        \REQUIRE User's query $Q$, Max iterations $T_{\max}$, LLM prompts $P^{\text{IG}}$ and $P_0^{\text{PP}}$.   
        \ENSURE Best found decision variables $\mathbf{y}_{\text{best}}$ and satisfaction score $S_{\text{best}}$.
        \STATE \textit{Initialization:} $\bar{\mathbf{y}}_{hist} \leftarrow \text{GetHistoricalAverage}(\mathcal{D})$, $t \leftarrow 0$, $S_{-1} \leftarrow -\infty$, $S_0 \leftarrow 0$. $\mathbf{y}_{\text{best}} \leftarrow \bar{\mathbf{y}}_{hist}$, $S_{\text{best}} \leftarrow S_0$.
        \STATE \textbf{Problem Matcher:} Determine an area-specific agent.
        \STATE \textbf{Indicator Generator:} Generate query-relevant objective function and code: $f(\mathbf{y};\mathbf{w})=\textbf{LLM}(Q;\mathcal{D}, P^{\text{IG}})$ .
        % \STATE \hspace{10pt} 2. Objective $f(\mathbf{y};\mathbf{w})$ code generation.
        \STATE \textbf{Code Safeguard:} Check the code generated by LLMs.
        \WHILE{$t < T_{\max}$ \textbf{and} $S_t > S_{t-1}$}
        \STATE $t \leftarrow t+1$
       \STATE \textbf{{a.}\ Problem Tailor:}
        \STATE \hspace{10pt} 1. Call empirical data $\mathcal{D}$ from a database.
        \STATE \hspace{10pt} 2. Output the set of names for remained decision variables $\hat{\mathbf{y}}_t$: $\hat{\mathbf{y}}_t=\textbf{LLM}(Q;\mathcal{D},P_t^{\text{PP}}(\hat{\mathbf{y}}_{t-1},S_{t-1})).$
        \STATE \hspace{10pt} 3. Add variable constraints for $ \mathbf{y}'_t= \mathbf{y} \setminus \hat{\mathbf{y}}_t$: $\mathbf{y}'_t \leftarrow \bar{\mathbf{y}}'_t$.\\
        \STATE \textbf{b.\ Solve the reduced bi-objective MIP:}
        \begin{center}
            $\hat{\mathbf{y}}^*_t \leftarrow \underset{\hat{\mathbf{y}_t}} {\text{max}} (g(\hat{\mathbf{y}}_t,\bar{\mathbf{y}}_t^{\prime};\mathbf{w}), f(\hat{\mathbf{y}}_t,\bar{\mathbf{y}}_t^{\prime};\mathbf{w}))$
        \end{center}
        \STATE \hspace{10pt} 1. Update optimal solutions: $\mathbf{y}^*_t \leftarrow (\hat{\mathbf{y}}^*_t, \bar{\mathbf{y}}'_t)$.
        \STATE \hspace{10pt} 2. Update $S_t=\frac{f(\hat{\mathbf{y}}^*_t,\bar{\mathbf{y}}^{\prime};\mathbf{w})-f(\bar{\mathbf{y}}_{hist};\mathbf{w})}{f(\bar{\mathbf{y}}_{hist};\mathbf{w})}$. 
        % \STATE \hspace{10pt} 3. Update $P_{t+1}^{\text{PP}}(\mathbf{y}^{\prime}_{t},S_{t})$.
        % \STATE \hspace{10pt} \textbf{return} $\hat{\mathbf{y}}^*_{t+1},S_{t+1}$.
        \IF{$S_t > S_{\text{best}}$}
                \STATE $S_{\text{best}} \leftarrow S_t$, $\mathbf{y}_{\text{best}} \leftarrow \mathbf{y}^*_t$.
        \ENDIF
        \ENDWHILE
        \RETURN $ \mathbf{y}_{best},\ S_{best}$.
	\end{algorithmic}\label{pse:1} 
\end{algorithm}
The central challenge in urban operations management is bridging the gap between high-level, often ambiguous, natural language objectives and the vast, combinatorial solution space of mathematical optimization. A brute-force approach is computationally intractable. The core intuition of RideAgent is to avoid this by establishing a synergistic, iterative dialogue between a heuristic \textit{guide} and a rigorous \textit{solver}.
This concept, which we term \textit{Small-Sample Guided Optimization}, is the cornerstone of our framework. Instead of treating the problem monolithically, we decompose it into two specialized roles:

1. Heuristic Guide (LLM-powered Problem Tailor): This module acts like an experienced human operations manager. It leverages the LLM's reasoning capabilities to interpret the user's query $Q$ in the context of historical data $\mathcal{D}$. It does not attempt to find the globally optimal solution itself. Instead, its purpose is to provide intelligent heuristics—to identify the most promising sub-region of the vast solution space and advise the solver on where to focus its efforts. This is achieved by selecting a subset of decision variables $\hat{\mathbf{y}}$ for optimization and fixing the remaining variables to pre-determined values.

2. Rigorous Solver (RF-based Feature-Driven Model): This module is highly effective at finding a mathematically verifiable optimal solution within a well-defined and reasonably sized problem space. By accepting the advice from the Heuristic Guide, it can direct its powerful search capabilities to the most promising areas and bypass the need to explore low-quality regions of the solution space.

The iterative loop between these two components creates a powerful feedback mechanism. The Guide proposes a search direction; the Solver explores it and reports back its findings; the satisfaction score $S_t$ quantifies the quality of this result and informs the Guide's next suggestion. This process allows RideAgent to efficiently navigate an immense decision space by combining the contextual understanding and flexible reasoning of LLMs with the mathematical precision of traditional optimization solvers.
\begin{figure*}
    \centering
    \includegraphics[width=\textwidth]{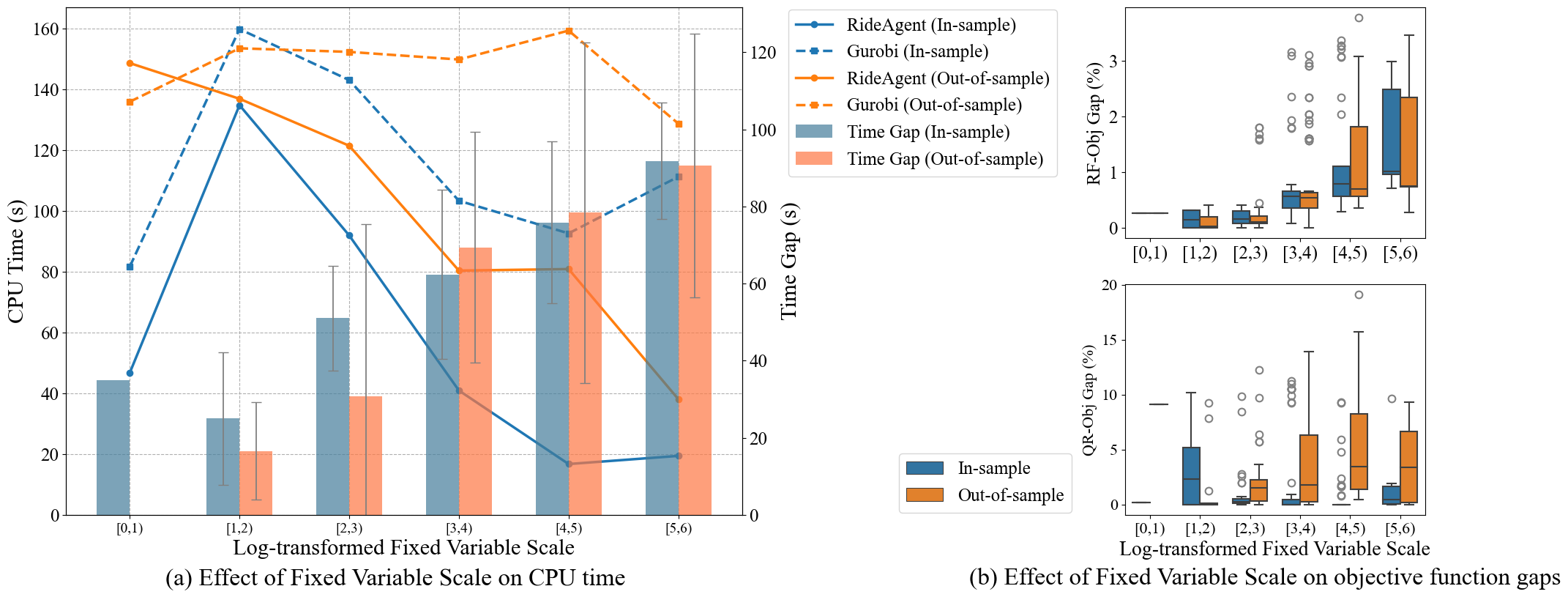}
    \captionsetup{name={Fig.},labelsep=period,justification=raggedright, singlelinecheck=false}
    \caption{Effect of \textit{Fixed Variable Scale} on linear objective function gaps}
    \label{fig:Linear Gaps}
\end{figure*}

\section{Case Study}
In this section, we conduct a case study using a real-world taxi dataset from New York City to validate the RideAgent framework. The numerical experiments are designed to demonstrate RideAgent's capabilities in translating diverse user objectives into actionable decisions and to evaluate the framework's computational efficiency.
%In this section, we conduct a case study using real-world taxi trip data. Numerical experiments demonstrate the multiple capabilities of RideAgent. %, including accuracy in generating indicators and high efficiency in solving problems. 
%Multiple experimental results reveal the great potential of intelligent agents in practical operation scenarios.
\subsection{Dataset Description }
Experiments are conducted using New York City yellow taxi trip records in 2016 ~\cite{nyctripdataset}. The trip records of New York yellow taxis include information such as trip start/end times and locations. Since the dataset lacks information on real-time taxi counts and historical dispatch decisions, we simulate this data using a flow-based model where taxis completing a trip in a zone become available supply for the next time step. Weather data for New York City in 2016, including features like temperature and dew point, are sourced from Kaggle ~\cite{weatherdataset}. The day of the week was also included as an exogenous feature.
The inconvenience cost, $\hat{w}_{ij}$, is set to \$0.5 per kilometer. The operator's fixed share, $\theta$, is 0.2 and the online booking fee $b_j= \$5$ per trip. Assume that electric taxis have three discrete power levels: low, medium, and high, corresponding to SOC of 0, 1, and 2 respectively.
For spatial analysis, all trip origins and destinations are clustered into 50 zones based on their geographical locations. The morning peak time is defined as 8:00 a.m. to 8:30 a.m. 
%In order to alleviate the mismatch between taxi supply and demand during the morning peak, the taxi pre-allocation operation and pricing decision are executed before the morning peak time. 
% According to the historical trip records and the simulation results of regional taxi inventory, it can be found that the taxi demand in the morning peak time in 8 areas is not fully met by the regional taxi inventory, while the remaining 42 areas have idle taxis. They are defined as demand and supply areas respectively.
Our analysis of historical data reveals a consistent supply-demand imbalance: 8 zones regularly exhibit unmet demand (defined as ``demand areas''), while the remaining 42 have surplus taxis (defined as ``supply areas'').
The RF model, which predicts operational profit, is trained on the morning peak operation data and exogenous features of 366 days in 2016. In the RF, the number of trees is 200 and the maximum depth is limited to 150. 30\% of the data is used as a test set. The RF achieves an R-squared value of 93.4\% and 60.7\% on the training set and test set, respectively.
To create a small sample of historical optimal decisions, we solve for the profit-maximizing (obj \ref{obj-rf}) allocation and pricing optimal decisions for 14 randomly selected days. These optimal decisions are provided as input to the agent as small sample optimal decision data.

All experiments are conducted on a macOS system equipped with an Apple M1 Pro CPU and 16GB of RAM. 
With GPT-4o mini as the core LLM, The agent framework is implemented in Python leveraging Langchain, a standard framework for creating applications powered by LLMs.
\subsection{Experimental Settings}
% We conduct three sets of experiments to evaluate the accuracy of indicator generation and the efficiency of model solving. 
% In all experiments, we adopt the unit-test methodology, common in software engineering and increasingly used for evaluating large agentic systems, to assess individual components through targeted test queries~\cite{li2023large}.
We evaluate the accuracy of objective generation and the efficiency of model solving using a unit-testing methodology, inspired by software engineering~\cite{li2023large}.
To ensure a comprehensive evaluation, we design a series of 18 queries, 15 of which can be formulated as linear objective functions, and the remaining three need to be formulated as nonlinear objective functions. The relevant objective function of each query is accompanied by a correct answer verified by human annotators. % (Further details in Appendix \ref{isa}).
To account for the inherent stochasticity of LLM outputs~\cite{casper2023open}, each query is executed 10 times. The reported accuracy and efficiency metrics are the average over these 10 runs. The underlying MIP models are solved using Gurobi 10.0.

\subsection{Experiment Results}
\subsubsection{Tests on the Accuracy of Objective Function Generation}
\begin{figure*}
\centering
\includegraphics[width=\textwidth]{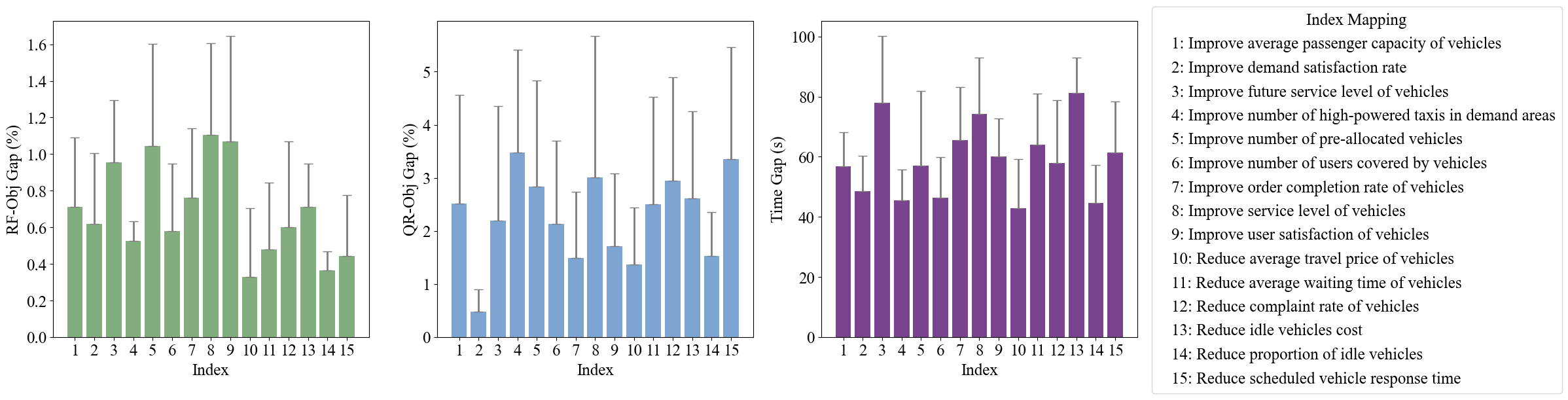} 
\captionsetup{name={Fig.},labelsep=period,justification=raggedright, singlelinecheck=false}
\caption{Agent performance on different user queries}
\label{index performance}
\end{figure*}
% To demonstrate the validity of the generated objective function, we execute accuracy tests by assessing the similarity between the objective functions generated by RideAgent and human-annotated standard objective functions.
% Sometimes the code of two functions may diverge, but their underlying mathematical significance could be the same. To address this problem, we introduce two metrics based on the Jaro-Winkler distance algorithm: \textit{text similarity} and \textit{result similarity}. These metrics quantify the similarity between the ``generated objective function" and the ``standard objective function" concerning their codes and mathematical essence, respectively (see Appendix \ref{essence} for details).
We first evaluate the ability of the \textbf{\textit{Indicator Generator}} to accurately translate a user's natural language query into a valid objective function. The assessment is performed by comparing the LLM-generated function against a human-annotated, ground-truth objective function. A key challenge is that functionally equivalent objectives can be expressed with syntactically different code. To account for this, we employ two distinct metrics based on the Jaro-Winkler distance: \textbf{\textit{text similarity}} and \textbf{\textit{result similarity}}. These metrics quantify the similarity between the ``generated objective function" and the ``standard objective function" concerning their codes and mathematical essence, respectively. The experiments are conducted in both an \textbf{\textit{in-sample} }setting (where the model is tested on a query it has seen as an example in its prompt) and an \textbf{\textit{out-of-sample}} setting (where the model is tested on a novel query).
% \begingroup
% \setlength\abovedisplayskip{0pt}
% \setlength\belowdisplayskip{0pt}
% \setlength\abovedisplayshortskip{0pt}
% \setlength\belowdisplayshortskip{0pt}
% \setlength\abovecaptionskip{5pt}
% \setlength\belowcaptionskip{0pt}
% \setlength{\tabcolsep}{1mm}
% \begin{table}[htbp]
% \centering
% \renewcommand{\arraystretch}{0.7}
% \caption{Accuracy Test Results for Linear Objectives}\label{table:relevance}
% \begin{threeparttable} 
% \begin{tablenotes}    
% \footnotesize              
% \item[1] The term ``Prompts" refers to the number of standard objective functions used in each experiment. 
% \end{tablenotes}   
% \begin{adjustbox}{max width=\columnwidth} % Adjust table width
% {\fontsize{10}{12}\selectfont\begin{tabular}{@{}cccccc@{}}
% \toprule
% {Prompts\tnote{1}} & \multicolumn{2}{c}{{In-sample test}} & \multicolumn{2}{c}{{Out-of-sample test}} \\ \cmidrule(l){2-5} 
%  & {Result similarity} & {Text similarity} & {Result similarity} & {Text similarity} \\ \midrule
% 0 & —— & —— & 0.72 & 0.86 \\
% 5 & 0.85 & 0.94 & 0.79 & 0.90 \\
% 10 & 0.89 & 0.95 & 0.81 & 0.91 \\
% 15 & 0.93 & 0.96 & —— & —— \\ \bottomrule
% \end{tabular}    
% }
% \end{adjustbox}

% \end{threeparttable}
% \end{table}

% \endgroup

\begingroup
\setlength\abovedisplayskip{0pt}
\setlength\belowdisplayskip{0pt}
\setlength\abovedisplayshortskip{0pt}
\setlength\belowdisplayshortskip{0pt}
\setlength\abovecaptionskip{5pt}
\setlength\belowcaptionskip{0pt}
\setlength{\tabcolsep}{4mm} % <-- 修改：加宽列间距，原为 1mm
\begin{table}[htbp]
\centering
\renewcommand{\arraystretch}{0.7}
\caption{Accuracy Test Results for Linear Objectives}\label{table:relevance}
\begin{threeparttable}
\begin{tablenotes}
\small % <-- 修改：加大标注字号，原为 \footnotesize
\item[1] The term ``Prompts" refers to the number of query-function pairs provided to the LLM in each experiment.
\end{tablenotes}
\begin{adjustbox}{max width=\columnwidth} % Adjust table width
{\fontsize{9}{12}\selectfont
\begin{tabular}{@{}cccccc@{}}
\toprule
{Prompts\tnote{1}} & \multicolumn{2}{c}{{In-sample}} & \multicolumn{2}{c}{{Out-of-sample}} \\ \cmidrule(l){2-5}
 & \begin{tabular}[c]{@{}c@{}}Result \\ similarity\end{tabular} & \begin{tabular}[c]{@{}c@{}}Text \\ similarity\end{tabular} & \begin{tabular}[c]{@{}c@{}}Result \\ similarity\end{tabular} & \begin{tabular}[c]{@{}c@{}}Text \\ similarity\end{tabular} \\ \midrule
0 & —— & —— & 0.72 & 0.86 \\
5 & 0.85 & 0.94 & 0.79 & 0.90 \\
10 & 0.89 & 0.95 & 0.81 & 0.91 \\
15 & 0.93 & 0.96 & —— & —— \\ \bottomrule
\end{tabular}
}
\end{adjustbox}
\end{threeparttable}
\end{table}
\endgroup

\begingroup
\setlength\abovedisplayskip{0pt}
\setlength\belowdisplayskip{0pt}
\setlength\abovedisplayshortskip{5pt}
\setlength\belowdisplayshortskip{5pt}
\setlength\abovecaptionskip{5pt}
\setlength\belowcaptionskip{0pt}
\setlength{\tabcolsep}{4mm} % <-- 修改：加宽列间距，原为 1mm
\begin{table}[htbp]
\centering
\renewcommand{\arraystretch}{0.7}
\caption{Accuracy Test Results for Nonlinear Objectives}
\begin{threeparttable}
% \begin{tablenotes}
% \footnotesize
% \item[1] The term ``Prompts" refers to the number of Q&A instructions used in each experiment.
% \end{tablenotes}
\begin{adjustbox}{max width=\columnwidth} % Adjust table width
{\fontsize{9}{12}\selectfont
\begin{tabular}{@{}cccccc@{}}
\toprule
{Prompts\tnote{1}} & \multicolumn{2}{c}{{\textbf{In-sample}}} & \multicolumn{2}{c}{{\textbf{Out-of-sample}}} \\ \cmidrule(l){2-5}
 & \begin{tabular}[c]{@{}c@{}}Result \\ similarity\end{tabular} & \begin{tabular}[c]{@{}c@{}}Text \\ similarity\end{tabular} & \begin{tabular}[c]{@{}c@{}}Result \\ similarity\end{tabular} & \begin{tabular}[c]{@{}c@{}}Text \\ similarity\end{tabular} \\ \midrule
0 & —— & —— & 0.71 & 0.86 \\
1 & 0.81 & 0.94 & 0.75 & 0.85 \\
2 & 0.83 & 0.95 & 0.78 & 0.84 \\
3 & 0.83 & 0.96 & —— & —— \\ \bottomrule
\end{tabular}
}
\end{adjustbox}
\end{threeparttable}
\label{table:nonlrelevance}
\end{table}
\endgroup

The results presented in Tables \ref{table:relevance} and \ref{table:nonlrelevance} lead to two key findings. 
First, RideAgent demonstrates strong \textbf{\textit{zero-shot capability}}. Even with zero examples provided in the prompt, it achieves a high text similarity of 86\% for both linear and nonlinear objectives, which indicates a robust intrinsic ability to interpret user requests. 
Second, the performance consistently improves as more examples are included in the prompt, which demonstrates effective \textbf{\textit{few-shot learning}}. For instance, in the linear out-of-sample test, result similarity increases from 0.72 to 0.81 as the number of prompts grows from 0 to 10. 
These results confirm that RideAgent can leverage the generalization capabilities of LLMs to reliably formulate objectives that align with user intent. Furthermore, they demonstrate that incorporating a small amount of domain-specific expert knowledge via few-shot prompting is an effective strategy to mitigate LLM hallucinations and enhance generation accuracy.

% Tables \ref{table:relevance} and \ref{table:nonlrelevance} show that RideAgent performs well in the out-of-sample test where no standard answer is given in the prompt ($0$ prompt), and the text similarity of both linear and non-linear indicator functions reaches 86\%.
% As the number of ``Prompts" increases, it can be observed that the similarity between the generated objective and the standard objective functions gradually increases in both in-sample and out-of-sample tests. The results indicate that RideAgent utilizes the generalization ability of LLMs to formulate objectives that are closely related to user queries. The results also show that the hallucination problem when LLMs are applied to objective function generation can be alleviated by adding expert knowledge to the prompts.
%Based on effective query-revelant objectives, we next explore whether RideAgent has advantages in reformulating and solving scaled-down mathematical programs, contrasting it with the full-scale optimization model.
% \begin{figure*}
%     \centering
%     \includegraphics[width=\textwidth]{Figure/Plot_Time_Performance_Advantage.png}
%     \caption{Reduced Solving Time by Agent}
%     \label{fig:time_gap}
% \end{figure*}
\begin{table*}[htbp]
\caption{Results of the Efficiency Test for Linear Objectives}\label{table:acc-linear}
\resizebox{\textwidth}{!}{%
{\fontsize{9}{12}\selectfont\begin{tabular}{@{}ccccccccc@{}}
\toprule
\multirow{4}{*}{\textbf{\begin{tabular}[c]{@{}c@{}}Fixed Variable \\ Scale\end{tabular}}} & \multicolumn{4}{c}{\textbf{In-sample}}   & \multicolumn{4}{c}{\textbf{Out-of-sample}}   \\ \cmidrule(l){2-9} 
& \multicolumn{3}{c}{RideAgent} & FULL  & \multicolumn{3}{c}{RideAgent} & FULL                                                                  \\ \cmidrule(l){2-9} 
& \multicolumn{2}{c}{Optimization GAP(\%)}      & \multirow{2}{*}{\begin{tabular}[c]{@{}c@{}}CPU time\\ (s)\end{tabular}} & \multirow{2}{*}{\begin{tabular}[c]{@{}c@{}}CPU time\\ (s)\end{tabular}} 
& \multicolumn{2}{c}{Optimization GAP(\%)}       & \multirow{2}{*}{\begin{tabular}[c]{@{}c@{}}CPU time\\ (s)\end{tabular}} & \multirow{2}{*}{\begin{tabular}[c]{@{}c@{}}CPU time\\ (s)\end{tabular}} 
\\ \cmidrule(lr){2-3} \cmidrule(lr){6-7}
 & \begin{tabular}[c]{@{}c@{}}RF-Obj gap\end{tabular} & \begin{tabular}[c]{@{}c@{}}QR-Obj Gap\end{tabular} & &
 & \begin{tabular}[c]{@{}c@{}}RF-Obj Gap\end{tabular} & \begin{tabular}[c]{@{}c@{}}QR-Obj Gap\end{tabular} &                 \\ \midrule
{[}0-50{]}& 0.41 & 1.73 & 65.19  & 120.07 & 0.57& 2.68  & 99.76  & 150.81  \\
{[}50-100{]}& 1.23 & 1.31 & 15.81 & 96.16 & 1.10  & 5.78 & 102.02& 159.34 \\
{[}100-150{]}& 1.17 & 0.05 & 16.69 & 85.75  & 1.20 & 5.78& 66.51 & 156.51    \\
{[}150-250{]} & 1.61  & 2.07 & 19.42 & 111.21& 1.47& 3.82& 38.08& 128.68   \\ \bottomrule
\end{tabular}
}}

\end{table*}
\begin{table*}[htbp]
\caption{Results of the Efficiency Test for Nonlinear Objectives}\label{table:acc-nonlinear}
\resizebox{\textwidth}{!}{%
{\fontsize{9}{12}\selectfont\begin{tabular}{@{}cllllllllllll@{}}
\toprule
\multirow{3}{*}{\textbf{\begin{tabular}[c]{@{}c@{}}Fixed Variable \\ Scale\end{tabular}}} & \multicolumn{4}{c}{Dispatching efficiency} & \multicolumn{4}{c}{Market share} & \multicolumn{4}{c}{Supply-demand matching degree} \\ \cmidrule(l){2-13} 
 & \multicolumn{2}{c}{CPU time(s)} & \multicolumn{1}{c}{\multirow{2}{*}{\begin{tabular}[c]{@{}c@{}}RF-Obj\\ Gap(\%)\end{tabular}}} & \multicolumn{1}{c}{\multirow{2}{*}{\begin{tabular}[c]{@{}c@{}}QR-Obj\\ Gap(\%)\end{tabular}}} & \multicolumn{2}{c}{CPU time(s)} & \multicolumn{1}{c}{\multirow{2}{*}{\begin{tabular}[c]{@{}c@{}}RF-Obj\\ Gap(\%)\end{tabular}}} & \multicolumn{1}{c}{\multirow{2}{*}{\begin{tabular}[c]{@{}c@{}}QR-Obj\\ Gap(\%)\end{tabular}}} & \multicolumn{2}{c}{CPU time(s)} & \multicolumn{1}{c}{\multirow{2}{*}{\begin{tabular}[c]{@{}c@{}}RF-Obj\\ Gap(\%)\end{tabular}}} & \multicolumn{1}{c}{\multirow{2}{*}{\begin{tabular}[c]{@{}c@{}}QR-Obj\\ Gap(\%)\end{tabular}}} \\ \cmidrule(lr){2-3} \cmidrule(lr){6-7} \cmidrule(lr){10-11}
& \multicolumn{1}{c}{Agent} & \multicolumn{1}{c}{FULL} & \multicolumn{1}{c}{} & \multicolumn{1}{c}{} & \multicolumn{1}{c}{Agent} & \multicolumn{1}{c}{FULL} & \multicolumn{1}{c}{} & \multicolumn{1}{c}{} & \multicolumn{1}{c}{Agent} & \multicolumn{1}{c}{FULL} & \multicolumn{1}{c}{} & \multicolumn{1}{c}{} \\ \midrule
10\% & 133.93 & 221.18 & 1.08 & 13.90 & 209.67 & 261.41 & 0.10 & 4.10 & 35.29 & 162.53 & 2.42 & 2.82 \\
20\% & 31.00 & 225.06 & 2.60 & 19.54 & 11.45 & 241.29 & 3.12 & 3.88 & 98.78 & 158.40 & 0.20 & 2.22 \\
30\% & 160.31 & 223.75 & 2.84 & 14.82 & 185.17 & 236.24 & 0.36 & 13.52 & 88.14 & 156.69 & 0.88 & 5.52 \\ \bottomrule
\end{tabular}
}}
\end{table*}

\subsubsection{Tests on the Efficiency of Model Solving}
In this section, we evaluate the core trade-off of the RideAgent framework: its ability to reduce computational time while maintaining near-optimal solution quality.
\paragraph{Motivation and Benchmark}
While existing LLM-assisted OR tools (\textit{OptiMUS}~\cite{ahmaditeshnizi2023optimus}, \textit{ORPO}~\cite{yang2023large}, \textit{OptiGuide}~\cite{li2023large}) have advanced human-computer interaction, they are not designed for the specific challenge of our case study, which requires integrating rich historical data to heuristically guide a large-scale feature-driven optimization model. To provide a rigorous benchmark, we compare RideAgent's performance against a \textbf{baseline model}, which we term the ``FULL'' model. This is the complete, feature-driven MIP model solved without any heuristic variable fixing from the Problem Tailor.

\paragraph{Evaluation Metrics}
\begin{itemize}
    \item \textbf{CPU Time}: The time required to solve the optimization problem.
    \item \textbf{Time Gap}: The reduction in computational time achieved by RideAgent relative to the FULL model. It is calculated as $(\text{CPU Time}_{\text{FULL}} - \text{CPU Time}_{\text{RideAgent}})$.
    \item \textbf{RF-Obj Gap}: The percentage deviation of RideAgent's primary objective (operational profit) from the FULL model's optimal profit.
    \item \textbf{QR-Obj Gap}: The percentage deviation of RideAgent's secondary, query-relevant objective from that of the FULL model.
\end{itemize}
Together, these metrics allow us to quantify the trade-off between computational efficiency (CPU Time) and solution optimality (the two objective gaps). For all tests, the number of few-shot examples (\textit{Prompts}) is fixed at 8. Define the metric \textit{Fixed Variable Scale} as the number of decision variables heuristically fixed by RideAgent. 

\paragraph{Performance on Linear Query-Relevant Objectives}
We first analyze the results for linear query-relevant objectives. Figure \ref{fig:Linear Gaps} illustrates the overall effect of the \textit{Fixed Variable Scale} while Table \ref{table:acc-linear} provides the aggregated numerical data.
Figure \ref{fig:Linear Gaps} provides a detailed analysis using a log-transformed x-axis to better visualize trends across different scales. Figure \ref{fig:Linear Gaps}(a) illustrates the clear and growing computational advantage of our approach. It reveals that RideAgent's solution time is consistently lower than the FULL model's, and crucially, this time advantage widens as more variables are fixed. This efficiency gain, however, introduces a trade-off in solution quality, which is detailed in Figure \ref{fig:Linear Gaps}(b). The box plots show that as the \textit{Fixed Variable Scale} increases, both the median objective gaps and their variance tend to increase. Furthermore, the in-sample tests (blue boxes) consistently show lower medians and tighter distributions than the out-of-sample tests (orange boxes). This finding demonstrates that providing relevant examples improves the reliability of the agent's heuristic guidance. Table \ref{table:acc-linear} quantifies the increasing gap trend with concrete numbers: the in-sample RF-Obj Gap climbs from 0.41\% at the lowest scale to 1.61\% at the highest. Simultaneously, it also confirms that even at its peak, this primary objective gap remains at an acceptably low level.

To provide a more granular view, Figure \ref{index performance} breaks down the performance for each of the 15 specific linear user queries. The RF-Obj Gap (left panel) remains low across all queries (typically below 1.2\%), and the Time Gap (right panel) is consistently large. The QR-Obj Gap (middle panel) exhibits more variability. Notably, for queries like \#4 (``Improve number of high-powered taxis") and \#12 (``Reduce complaint rate"), the QR-Obj Gap is extremely low. This suggests that the historical data provided strong guidance for these particular objectives.

\begin{table*}[htbp]
\centering
\caption{Comparison Results with Cutting Methods}
\label{table:cuts}
\resizebox{\textwidth}{!}{%
\begin{tabular}{@{}ccccccc@{}}
\toprule
\multirow{2}{*}{\textbf{Cuts Name}} & \multicolumn{3}{c}{\textbf{In-sample}} & \multicolumn{3}{c}{\textbf{Out-of-sample}} \\ \cmidrule(l){2-7} 
 & Time Gap & RF-Obj Gap & QR-Obj Gap & Time Gap & RF-Obj Gap & QR-Obj Gap \\ \midrule
CliqueCuts~\cite{atamturk2000conflict} & 39.46s(30.49\%) & 0.80\% & 5.13\% & 50.75s(42.58\%) & 1.08\% & 5.30\% \\
CoverCuts~\cite{ceria1998cutting} & 43.18s(33.09\%) & 0.80\% & 5.13\% & 55.44s(46.75\%) & 1.08\% & 5.09\% \\
GomoryCuts~\cite{gomory2010outline} & 53.00s(43.01\%) & 0.80\% & 5.26\% & 62.33s(52.48\%) & 1.08\% & 4.97\% \\
GUBCoverCuts~\cite{wolsey1990valid} & 46.45s(36.80\%) & 0.80\% & 5.18\% & 59.35s(49.94\%) & 1.08\% & 5.02\% \\
MIRCuts~\cite{marchand2001aggregation} & 46.56s(36.46\%) & 0.80\% & 5.13\% & 61.31s(51.65\%) & 1.08\% & 5.27\% \\ \midrule
\textbf{Total average} & 45.73s(35.97\%) & 0.80\% & 5.16\% & 57.83s(48.68\%) & 1.08\% & 5.13\% \\ \bottomrule
\end{tabular}%
}
\end{table*}

\begin{figure*}
\centering
\includegraphics[width=\textwidth]{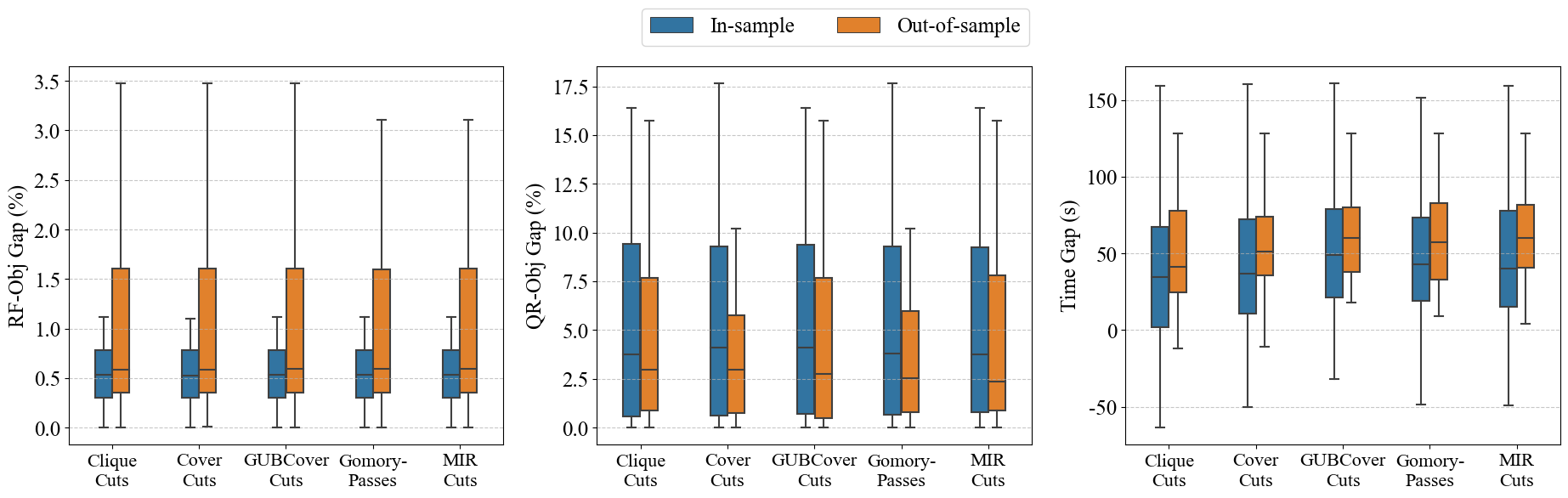} 
\captionsetup{name={Fig.},labelsep=period,justification=raggedright, singlelinecheck=false}
\caption{Objective function gaps and time advantage of RideAgent compared to different cutting methods}
\label{Cut gaps}
\end{figure*}

\paragraph{Performance on Nonlinear Query-Relevant Objectives}
The results for more computationally intensive nonlinear queries are summarized in Table \ref{table:acc-nonlinear}. For these tests, the \textit{Fixed Variable Scale} are set to specific percentages (10\%, 20\%, 30\%) of the 1032 total decision variables.
RideAgent's advantage is even more pronounced in this context. The average CPU time is reduced by a remarkable 50.04\%, while the RF-Obj Gap remains exceptionally low at an average of just 1.51\% across all tests. However, the QR-Obj Gap is significantly higher than in the linear case, and it increases substantially as more variables are fixed. This result provides a crucial insight into the framework's mechanics. The historical guidance data is optimized for the primary profit objective (RF-Obj). When a user introduces a conflicting, nonlinear objective, this guidance is less relevant. As a result, the deviation from the optimal solution for that specific query is larger. This highlights that the effectiveness of the heuristic depends on the alignment between the historical guidance data and the new user query.

\subsubsection{Comparison with Cutting Planes Methods}
% To further benchmark the effectiveness of our LLM-guided heuristic, we compare it with five commonly used cuts: CliqueCuts, CoverCuts, GomoryPassesCuts, GUBCoverCuts, MIRCuts. 
% For each cut type, we establish a benchmark by solving the FULL model with only that specific cut disabled. We then compare the performance (optimality gaps and CPU time) of this ``cuts-off" benchmark against the performance of RideAgent. This allows for a direct comparison between two distinct strategies for accelerating the solution process: a classical OR technique versus our LLM-guided heuristic. To ensure a fair comparison, only linear query-relevant objectives were used, and the number of prompts is fixed at 8.
We benchmark our LLM-guided heuristic against five standard cutting planes (CliqueCuts, CoverCuts, GomoryPassesCuts, GUBCoverCuts, MIRCuts) on linear objectives. The performance of RideAgent is compared against solving the FULL model with only one of these cut types disabled at a time, providing a direct comparison of solution acceleration strategies.

Figure \ref{Cut gaps} presents a visual comparison of the performance trade-offs, while Table \ref{table:cuts} provides the corresponding mean values. The left panel of Figure \ref{Cut gaps} shows that the RF-Obj Gap for RideAgent is consistently low, with a tight distribution, especially for the in-sample (blue) tests. In contrast, the middle and right panels display a more complex, counter-intuitive pattern. 
Specifically, in-sample tests result in a surprisingly higher median and variance for the QR-Obj Gap (middle panel), alongside a smaller average Time Gap (right panel). Both phenomena stem from the aggressive nature of the heuristic generated by the LLM for high-confidence, in-sample queries. This strategy retains profit-critical variables as active and creates a reduced model that is ideal for the primary profit objective but challenging for any conflicting secondary goal. This has two direct consequences: first, in the lexicographical second step, the solver is forced into a compromise further from the QR-Obj optimum, thus increasing its gap. Second, this same difficult search process increases RideAgent's solution time, which shrinks the average Time Gap. The high variance in both metrics is also explained by this ``high-risk, high-reward" heuristic. If the heuristic aligns with the user's query, the problem is solved almost instantly, which yields a large Time Gap. If the two conflict, the solver struggles with the difficult second-stage problem, which leads to a very long solution time and thus a small or even negative Time Gap. Conversely, out-of-sample queries yield more conservative heuristics and lead to lower variance.

Despite these dynamics, the Time Gap box plots are almost entirely in positive territory, which indicates that RideAgent consistently outperforms disabling a standard cut. Table \ref{table:cuts} quantifies these visual trends with precise numbers. On average, the LLM-guided strategy achieves a time saving of 51.78 seconds (42.32\%) while incurring only minor objective gaps (an average RF-Obj Gap of 0.94\% and QR-Obj Gap of 5.15\%). These results demonstrate a favorable balance between computational efficiency and solution quality.

\section{Conclusion}
This paper introduces \textbf{RideAgent}, a novel agent framework designed to bridge the persistent gap between ambiguous human objectives and the precise requirements of mathematical optimization. Our core contribution is a new paradigm we term \textit{Small-Sample Guided Optimization,} where an LLM acts as a heuristic guide for a feature-driven MIP model. The LLM interprets user intent and, guided by a small set of historical optimal examples, prunes the vast decision space so that the rigorous MIP solver can focus its computational power on the most promising regions.

Our case study provides robust empirical validation for this paradigm. The results demonstrate that RideAgent is not only conceptually novel but also practically effective.
First, its \textit{Indicator Generator} reliably translates natural language into valid objectives, and it achieves up to 86\% text similarity in zero-shot settings.
Second, the framework demonstrates remarkable efficiency. It reduces CPU time by over 80\% for complex nonlinear problems while keeping the primary objective gap below 1.5\%.
Finally, its heuristic strategy proves more effective than standard MIP acceleration techniques, and it outperforms five established cutting plane methods by an average of 52 seconds per solve.

Beyond the immediate application of taxi operations, we believe the \textit{Small-Sample Guided Optimization} framework offers a generalizable and powerful approach for a wide range of complex decision-making problems in domains such as supply chain management, logistics, and resource scheduling. By creating an effective bridge between human intuition and mathematical optimization, this paradigm promises to democratize the power of operations research and make it more accessible and responsive to the dynamic needs of practitioners.

Future research could proceed in several exciting directions. Enhancing the underlying predictive models would further improve the quality of the primary objective. Investigating more sophisticated dialogue protocols between the LLM and the solver could unlock further efficiencies. Finally, deploying and testing the RideAgent framework in other real-world operational environments will be a crucial next step in validating its broader applicability and impact.
\bibliography{ref}
\bibliographystyle{plainnat}
% \clearpage
\appendix
% \pagenumbering{arabic}
% \setcounter{table}{0}   
% \setcounter{figure}{0}
\section{Appendix}

\subsection{Objective Function Standard Answer}\label{isa}
In this section, we provide the ground-truth Q\&A instructions adopted in the prompt of RideAgent. 
Specifically, we provide accurate Q\&A information (human-labeled) in three specific scenarios: \emph{Operational}, \emph{Customer-Related}, and \emph{Regulatory}. These scenarios encompass the majority of operational inquiries proposed by taxi fleet operators. 
Within each scenario, there are typically 3 to 8 key QR-objs that are of utmost relevance to the majority of taxi fleet operators (as shown in Table \ref{app-tab:groundtruth}). 
For each QR-obj, we provide a ground-truth ``Answers" that is derived from human experience. This information is supplied in the form of a Gurobi objective code.

\subsection{Objective Function Similarity Index}\label{essence}
This section presents the mathematical representation of the Gurobi code created by RideAgent, which is utilized to quantify the ``Results Similarity" in the relevance test.
The Gurobi code is shown below, which represents the objective of maximizing the total number of accessible e-bikes and its related mathematical essence:

\begin{itemize}
    \item \textbf{Gurobi code}: \texttt{model.setObjective( gp.quicksum(model.getVarByName( f"cluster\{i\}\_cluster\{j\}\_\{k\}") for i in S for j in D for k in K), GRB.MAXIMIZE)}
    \item  \textbf{Mathematical essence}:  $\max\  \hat{u}_{0,0}+\hat{u}_{0,1}+\hat{u}_{0,2}+\underbrace{\quad\quad\quad \quad\cdots\cdots\quad\quad\quad\quad}_{\hat{u}_{j,k}\ for\  j\ in\ [4, 7, 12, 15, 37, 38, 44]\ for\ k\ in\ [0,1,2]}+\hat{u}_{49,0}+\hat{u}_{49,1}+\hat{u}_{49,2}$

    \item \textbf{Gurobi code}: \texttt{model.setObjective( gp.quicksum(model.getVarByName( f"cluster\{i\}\_cluster\{j\}\_\{k\}") for i in S for j in D for k in K if k >0), GRB.MAXIMIZE)}
    \item  \textbf{Mathematical essence}: $\max\  \hat{u}_{0,1}+\hat{u}_{0,2}+\underbrace{\quad\quad\quad \quad\cdots\cdots\quad\quad\quad\quad}_{\hat{u}_{j,k}\ for\  j\ in\ [4, 7, 12, 15, 37, 38, 44]\ for\ k\ in\ [0,1,2]}+\hat{u}_{49,1}+\hat{u}_{49,2}$

\end{itemize}
\subsection{Prompt}\label{prompt}
\begin{figure*}[htbp]
\centering
\includegraphics[width=\textwidth]{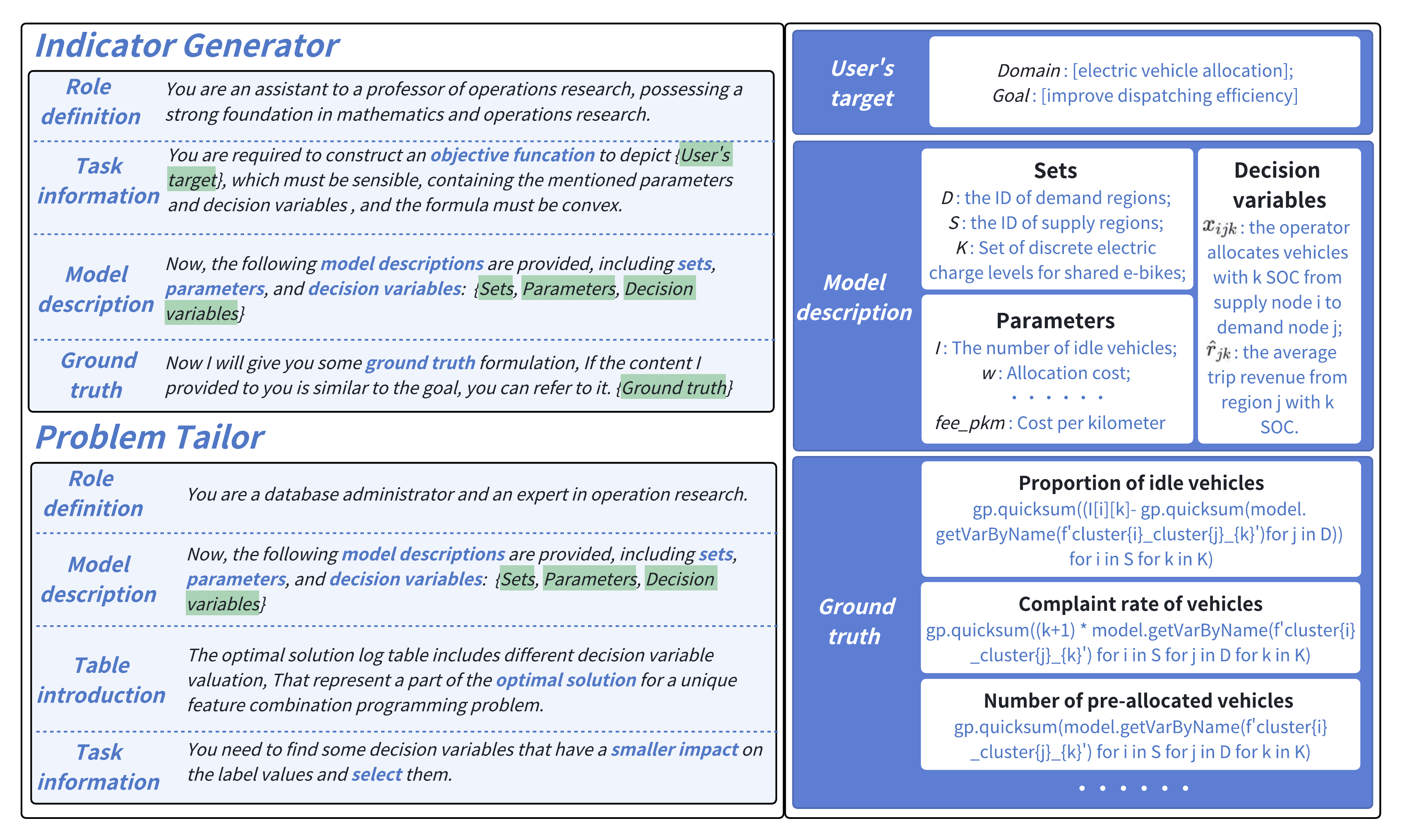} 
\captionsetup{name={Fig.},labelsep=period,justification=raggedright, singlelinecheck=false}
\caption{Brief introduction of prompt}
\label{Figure prompt}
\end{figure*}

Figure \ref{Figure prompt} shows the prompts input to the LLMs for two parts of the agent, \emph{Indicator Generator} and \emph{Problem Tailor}. These prompts contain role definitions, task descriptions, and valid problem information.

\begin{table*}[htbp]
\captionsetup{labelfont=bf, textfont=bf}
\caption{Ground-truth Objective Functions}\label{app-tab:groundtruth}
\centering
\renewcommand{\arraystretch}{2.5} % Increased spacing for better readability
\begin{adjustbox}{max width=\textwidth, center}
\Large % Increased font size
\begin{tabular}{ll}
\toprule
\textbf{Query Relevant Objective Function} & \textbf{Ground-truth Objective Function Code} \\
\midrule
\makecell[l]{Proportion of idle taxis} & \makecell[l]{$\mathtt{model.setObjective(gp.quicksum(S[i, k]\ -\ gp.quicksum(x[i, j, k]\ for\ j\ in\ D)}$\\$\mathtt{for\ i\ in\ S\ for\ k\ in\ K), GRB.MINIMIZE)}$} \\
\makecell[l]{Idle taxis cost} & \makecell[l]{$\mathtt{model.setObjective(gp.quicksum((k+1)\ *\ (S[i, k]\ -\ gp.quicksum(x[i, j, k]\ for\ j\ in\ D))}$\\ $\mathtt{for\ i\ in\ S\ for\ k\ in\ K), GRB.MINIMIZE)}$} \\
\makecell[l]{Number of high-powered taxis \\ in demand areas} & $\mathtt{model.setObjective(gp.quicksum(x[i, j, 2]\ for\ i\ in\ S\ for\ j\ in\ D), GRB.MAXIMIZE)}$ \\
\makecell[l]{Future service Level of taxis} & $\mathtt{model.setObjective(gp.quicksum(k\ *\ x[i, j, k]\ for\ i\ in\ S\ for\ j\ in\ D\ for\ k\ in\ K), GRB.MAXIMIZE)}$ \\
\makecell[l]{Scheduled taxi response time} & $\mathtt{model.setObjective(gp.quicksum(d[i, j]\ *\ x[i, j, k]\ for\ i\ in\ S\ for\ j\ in\ D\ for\ k\ in\ K), GRB.MINIMIZE)}$ \\
\makecell[l]{Dispatching efficiency of taxis} & \makecell[l]{$\mathtt{model.setObjective(gp.quicksum((u[j, k]\ -\ w[i, j])\ *\ x[i, j, k]}$\\ $\mathtt{for\ i\ in\ S\ for\ j\ in\ D\ for\ k\ in\ K), GRB.MAXIMIZE)}$} \\
\makecell[l]{Complaint rate of taxis} & $\mathtt{model.setObjective(gp.quicksum((k+1)\ *\ x[i, j, k]\ for\ i\ in\ S\ for\ j\ in\ D\ for\ k\ in\ K), GRB.MAXIMIZE)}$ \\
\makecell[l]{Service Level of taxis} & $\mathtt{model.setObjective(gp.quicksum((k+1)\ *\ x[i, j, k]\ for\ i\ in\ S\ for\ j\ in\ D\ for\ k\ in\ K), GRB.MAXIMIZE)}$ \\
\makecell[l]{Average travel price of taxis} & $\mathtt{model.setObjective(gp.quicksum(u[j, k]\ for\ j\ in\ D\ for\ k\ in\ K), GRB.MINIMIZE)}$ \\
\makecell[l]{Order completion rate of taxis} & $\mathtt{model.setObjective(gp.quicksum((k+1)\ *\ x[i, j, k]\ for\ i\ in\ S\ for\ j\ in\ D\ for\ k\ in\ K), GRB.MAXIMIZE)}$ \\
\makecell[l]{Average waiting time of taxis} & $\mathtt{model.setObjective(gp.quicksum(d[i, j]\ *\ x[i, j, k]\ for\ i\ in\ S\ for\ j\ in\ D\ for\ k\ in\ K), GRB.MINIMIZE)}$ \\
\makecell[l]{Supply-demand matching degree \\ of taxis} & \makecell[l]{$\mathtt{model.setObjective(gp.quicksum(S[i, k]\ -\ gp.quicksum(x[i, j, k]\ for\ j\ in\ D)\ for\ i\ in\ S\ for\ k\ in\ K)}$\\ $\mathtt{+\ gp.quicksum(gp.abs(demand\_avg[j]\ -\ inventory\_avg\ -}$\\ $\mathtt{gp.quicksum(x[i, j, k]\ for\ i\ in\ S\ for\ k\ in\ K))\ for\ j\ in\ D), GRB.MINIMIZE)}$} \\
\makecell[l]{Number of pre-allocated taxis} & $\mathtt{model.setObjective(gp.quicksum(x[i, j, k]\ for\ i\ in\ S\ for\ j\ in\ D\ for\ k\ in\ K), GRB.MAXIMIZE)}$ \\
\makecell[l]{Average passenger capacity \\ of taxis} & $\mathtt{model.setObjective(gp.quicksum((k+1)\ *\ x[i, j, k]\ for\ i\ in\ S\ for\ j\ in\ D\ for\ k\ in\ K), GRB.MAXIMIZE)}$ \\
\makecell[l]{Number of users covered \\ by taxis} & $\mathtt{model.setObjective(gp.quicksum((k+1)\ *\ x[i, j, k]\ for\ i\ in\ S\ for\ j\ in\ D\ for\ k\ in\ K), GRB.MAXIMIZE)}$ \\
\makecell[l]{User satisfaction of taxis} & $\mathtt{model.setObjective(gp.quicksum((k+1)\ *\ x[i, j, k]\ for\ i\ in\ S\ for\ j\ in\ D\ for\ k\ in\ K), GRB.MAXIMIZE)}$ \\
\makecell[l]{Demand satisfaction rate} & $\mathtt{model.setObjective(gp.quicksum((k+1)\ *\ x[i, j, k]\ for\ i\ in\ S\ for\ j\ in\ D\ for\ k\ in\ K), GRB.MAXIMIZE)}$ \\
\makecell[l]{Market share of taxis} & \makecell[l]{$\mathtt{model.setObjective(gp.quicksum(u[j, k]\ *\ x[i, j, k]\ for\ i\ in\ S\ for\ j\ in\ D\ for\ k\ in\ K), GRB.MAXIMIZE)}$} \\
\bottomrule
\end{tabular}
\end{adjustbox}
\end{table*}
\clearpage
% \section*{Acknowledgments}
% This should be a simple paragraph before the References to thank those individuals and institutions who have supported your work on this article.


\newpage



% \end{IEEEbiography}



% \vfill

\end{document}